\documentclass{amsart}
\usepackage{epsfig}
\usepackage{graphicx}
\usepackage{color}
\usepackage{amssymb}
\usepackage{esint}
\usepackage{tikz-cd}
\usetikzlibrary{matrix,arrows.meta}

\newcommand\R{\mathbb R}
\newcommand\C{\mathbb C}
\newcommand\N{\mathbb N}

\newcommand{\vertiii}[1]{{\left\vert\kern-0.25ex\left\vert\kern-0.25ex\left\vert #1
    \right\vert\kern-0.25ex\right\vert\kern-0.25ex\right\vert}}
\DeclareMathOperator\Tr{Tr}

\newtheorem{thm}{Theorem}

\newtheorem{lem}[thm]{Lemma}

\theoremstyle{definition}
\newtheorem{definition}[thm]{Definition}
\theoremstyle{remark}
\newtheorem{rem}{Remark}

\begin{document}
\title[On the eigenvalue problem for a bulk/surface ...]
{On the eigenvalue problem for a bulk/surface elliptic system}
\author{Enzo Vitillaro}
\address[E.~Vitillaro]
       {Dipartimento di Matematica e Informatica, Universit\`a di Perugia\\
       Via Vanvitelli,1 06123 Perugia ITALY}
\email{enzo.vitillaro@unipg.it}




\thanks{This work has been funded by the European Union - NextGenerationEU within the framework of PNRR  Mission 4 - Component 2 - Investment 1.1 under the Italian Ministry of University and Research (MUR) programme "PRIN 2022" - 2022BCFHN2 - Advanced theoretical aspects in PDEs and their applications - CUP: J53D23003700006}

\begin{abstract} The paper addresses the doubly elliptic eigenvalue problem
$$\begin{cases} -\Delta u=\lambda u \qquad &\text{in
$\Omega$,}\\
u=0 &\text{on $\Gamma_0$,}\\
-\Delta_\Gamma u +\partial_\nu u =\lambda u\qquad
&\text{on
$\Gamma_1$,}
\end{cases}
$$
where $\Omega$ is a bounded open subset of $\R^N$ ($N\ge 2$) with  a $C^1$ boundary  $\Gamma=\Gamma_0\cup\Gamma_1$, $\Gamma_0\cap\Gamma_1=\emptyset$, $\Gamma_1$ being nonempty and relatively open on $\Gamma$.
Moreover, $\mathcal{H}^{N-1}(\overline{\Gamma}_0\cap\overline{\Gamma}_1)=0$ and $\mathcal{H}^{N-1}(\Gamma_0)>0$.
We prove that $L^2(\Omega)\times L^2(\Gamma_1)$ admits a Hilbert basis constituted by eigenfunctions and  we describe the behavior of the eigenvalues. Moreover, when $\Gamma$ is at least $C^2$ and $\overline{\Gamma}_0\cap\overline{\Gamma}_1=\emptyset$, we give several qualitative properties of the eigenfunctions.
\end{abstract}

\maketitle

\section{Introduction and main results} \label{section 1}
\subsection{Presentation of the problem and literature overwiew}
We deal with the doubly elliptic eigenvalue problem
\begin{equation}\label{1}
\begin{cases} -\Delta u=\lambda u \qquad &\text{in
$\Omega$,}\\
u=0 &\text{on $\Gamma_0$,}\\
-\Delta_\Gamma u +\partial_\nu u =\lambda u\qquad
&\text{on
$\Gamma_1$,}
\end{cases}
\end{equation}
where $\Omega$ is a bounded open subset of $\R^N$ ($N\ge
2$) with $C^r$ boundary $\Gamma$ (see \cite{grisvard}), with $r=1,2,\ldots,\infty$. Hence, when nothing is said, $r=1$.
We also assume that $\Gamma=\Gamma_0\cup\Gamma_1$, $\Gamma_0\cap\Gamma_1=\emptyset$,
$\Gamma_1$ being nonempty and relatively open on $\Gamma$ (or equivalently $\overline{\Gamma_0}=\Gamma_0$).
Denoting by
$\mathcal{H}^{N-1}$ the  Hausdorff  measure, we also assume that   $\mathcal{H}^{N-1}(\overline{\Gamma}_0\cap\overline{\Gamma}_1)=0$ and $\mathcal{H}^{N-1}(\Gamma_0)>0$. These
properties of $\Omega$, $\Gamma_0$ and $\Gamma_1$ will be assumed,
without further comments, throughout the paper.

In problem \eqref{1} $\lambda$ is a real or complex parameter and we respectively denote by $\Delta$ and $\Delta_\Gamma$ the Laplace operator in $\Omega$ and the Laplace--Beltrami operator  on $\Gamma$, while $\nu$ stands for the outward normal to $\Omega$. We shall look for eigenvalues and eigenfunctions of problem \eqref{1}, that is for  values of $\lambda$ for which \eqref{1} has a nontrivial (real or complex--valued) solution, i.e. an eigenfunction.

Problem \eqref{1} has been studied (as the particular case $K=0$,  $\alpha=1$ and $\gamma=\omega$ in problem (1.2))  in
\cite{KnopfLiu}, when $\Gamma_0=\emptyset$ and $\lambda, u$ are real. 
The study in  \cite{KnopfLiu} is motivated by  several papers on the Allen--Cahn equation subject to a dynamic boundary condition, see
\cite{CalatroniColli, ColliFukao, SprekelsWu}. Indeed finding a Hilbert basis of eigenfunctions allows to look for  solutions of the evolution problem by using a Faedo--Galerkin scheme.

We remark that, by introducing  an (inessential) positive parameter $\kappa$ in front of the Laplacian in \eqref{1}, and formally taking the limit as $\kappa\to\infty$, one gets from \eqref{1} the  Wentzell eigenvalue problem, studied in \cite{Dambrine,DuWangXia}, which is then related to \eqref{1}.

As to the author's knowledge  problem \eqref{1}  in the case $\Gamma_0\not=\emptyset$ has not yet been considered in the mathematical literature. The motivation for studying it originates from an evolution problem different than the one mentioned above.
 Indeed, it originates from the wave equation with hyperbolic  boundary conditions. It is the evolutionary boundary value problem
\begin{equation}\label{2}
\begin{cases} w_{tt}-\Delta w=0 \qquad &\text{in
$\R\times\Omega$,}\\
w=0 &\text{on $\R\times \Gamma_0$,}\\
w_{tt}-\Delta_\Gamma
w +\partial_\nu w=0\qquad &\text{on
$\R\times \Gamma_1$,}
\end{cases}
\end{equation}
where $\Omega$, $\Gamma_0$ and $\Gamma_1$ are as above, $w=w(t,x)$, $t\in\R$, $x\in\overline{\Omega}$,  $\Delta=\Delta_x$ and  $\Delta_\Gamma$ denote the Laplacian and Laplace--Beltrami operators with respect to the space variable.

One easily see that solutions of \eqref{2} enjoy energy conservation, once a properly defined energy function is introduced.
So, while one cannot expect decay of solutions, it is of interest to look for standing wave solutions of \eqref{2}. They are solutions of the form
\begin{equation}\label{3}
w(t,x)=e^{i\omega t} u(x),\quad\omega \in\R\setminus \{0\},
\end{equation}
where $u$ is nontrivial and real--valued. The function $w$ defined in \eqref{3} solves, at least formally, problem \eqref{2} if and only if $u$ solves problem \eqref{1} with $\lambda=\omega^2>0$.

Hence, in the analysis of problem \eqref{2}, a deep understanding of the eigenvalue problem \eqref{1} would allow to find solutions by suing the Fourier method, provided one can find a complete system of eigenfunctions. Since in the analysis of problem (2) it is usefull to consider complex--valued solutions, in the sequel  we shall consider complex--valued functions everywhere. The due attention will be given to find real eigenfunctions, so a reader only interested to the real case can simply ignore conjugation everywhere.

Since our study of problem \eqref{1} is motivated by problem \eqref{2}, it is worth to give  a brief overview of the literature dealing with it. Indeed, problems with hyperbolic or, more generally, kinetic boundary conditions, arise in several physical applications.

 A one dimensional
model was studied by several authors
to describe transversal small oscillations of an elastic rod with a
tip mass on one endpoint, while the other one is pinched.
See \cite{andrewskuttlershillor,conradmorgul,darmvanhorssen,guoxu,morgulraoconrad}.

 A two dimensional model, 
\begin{footnote}{Although this model is two--dimensional, the title of the paper comes the model considered in \cite{KnopfLiu}, which is three dimensional, since the two models are formally strongly related.}\end{footnote}
introduced in \cite{goldsteingisele}, more closely motivates problem \eqref{2}.
We shall briefly describe it.
One considers a vibrating membrane of surface
density $\mu>0$, subject to a tension $T>0$, both taken constant and normalized  for simplicity.  Moreover, $w=w(t,x)$, $t\in\R$, $x\in\Omega\subset\R^2$,  denotes the vertical displacement from the
rest state. After a standard linear approximation, $w$
satisfies  the wave equation $w_{tt}-\Delta w=0$ in $\R\times\Omega$.
Now one supposes that a part $\Gamma_0$ of the
boundary  is pinched, while the other part $\Gamma_1$ carries a
 linear mass density $m>0$ and  it is subject to a linear
tension $\tau>0$, both taken constant for simplicity.

 A practical example of this situation is given by a
drumhead with a hole in the interior having  a thick border. This situation occurs
 in bass drums, as one can realize by looking at several pictures of them in the internet.
During the paper we shall constantly refer to this motivating  example as to \emph{the bass drum model}.

After a further linear approximation the boundary
condition thus reads as $mw_{tt}+\partial_\nu w-\tau
\Delta_{\Gamma}w=0$. In \cite{goldsteingisele} the case
$\Gamma_0=\emptyset$ and $\tau=0$ was considered, while here we
shall deal with the more realistic case $\Gamma_0\not=\emptyset$ and
$\tau>0$, with $\tau$ and $m$ normalized for simplicity.

We also would like to point out that, when $\Gamma_0=\emptyset$ and $\Omega=\R^N_+$, problem \eqref{2} also shows up in  Quantum Field Theory, see \cite{zahn}.

Several papers in the literature address the wave equation  with
kinetic boundary conditions. This fact is even more evident  if one
takes into account that, plugging the wave equation in \eqref{2} into the
boundary condition, we can rewrite it as  $-\Delta_\Gamma w +\partial_\nu
w+\Delta w=0$, that is a generalized Wentzell (also spelled Ventcel) boundary
condition.  We refer to  \cite{doroninlarkinsouza,FGGGR,MugnoloCosine,mugnolo2011,xiaoliang2,xiaoliang1}, and also
to the series of papers \cite{AMS,Dresda1,Dresda2,Dresda3} by the author.
All the quoted papers deal with well--posedness issues for variously (linearly or nonlinearly)  perturbed versions of \eqref{2}.

The stability issue for a damped (both internally and at the boundary) version of \eqref{2} was studied in
\cite{CavalcantiKM}, while a boundary damped version of it was subject of several papers, see \cite{heminna1,Heminna2,NicaiseLaoubi} and the more recent papers \cite{LLXiao1,LLXiao2,LiXiao}.

The analysis of the literature made above shows that the simple unperturbed problem \eqref{2},  which
  well--posedness is  rather standard, see for example \cite{Dresda1}, was  object of a detailed study only
  in the case $\Gamma_0=\emptyset$ and $\Omega=\R^N_+$, see \cite{zahn}.

  The aim of the present paper is to start such a study. In particular we shall analyze the eigenvalue problem \eqref{1}. Our aim is to show that most of the well--known classical results concerning the homogeneous Dirichlet problem for the Helmholtz equation, that is the eigenvalue problem
\begin{equation}\label{4}
-\Delta u=\lambda u \quad \text{in
$\Omega$,}\qquad\qquad u=0 \quad \text{on $\partial\Omega$,}
\end{equation}
continue to hold  for problem \eqref{1}, in a suitably modified form. The parallelism will go further than one can expect.
 To illustrate the last assertion we shall present our main results in the sequel. Preliminarily we are now going to introduce some basic notation.
\subsection{Functional spaces}In the paper we shall adopt the standard notation for (complex) Lebesgue and Sobolev spaces in $\Omega$ and on $\Gamma$, see \cite{adams} and \cite{grisvard}. We make the reader aware that Lebesgue spaces on $\Gamma$ are Lebesgue spaces with respect to the restriction of the Hausdorff measure $\mathcal{H}^{N-1}$ to measurable subsets of $\Gamma$. We shall  drop the notation
$d\mathcal{H}^{N-1}$ in boundary integrals, so for example $\int_{\Gamma_1}u= \int_{\Gamma_1}u\,d\mathcal{H}^{N-1}$. Moreover, we shall
identify $L^p(\Gamma_1)$, for $1\le p\le \infty$, with its isometric image in $L^p(\Gamma)$,
that is
\begin{equation}\label{5}
L^p(\Gamma_1)=\{u\in L^p(\Gamma):\,\, u=0\,\,\text{a.e. on
}\,\,\Gamma_0\}.
\end{equation}
Moreover, for $u\in H^1(\Omega)$,  we shall denote by $u_{|\Gamma}\in H^{1/2}(\Gamma)$ its trace on $\Gamma$.

We introduce the Hilbert space $H^0 = L^2(\Omega)\times L^2(\Gamma_1)$, endowed with the standard inner product
\begin{equation}\label{6}
  \Big((u_1,v_1),(u_2,v_2)\Big)_{H^0}=\int_\Omega u_1\overline{u_2}+\int_{\Gamma_1}v_1\overline{v_2}\quad\text{for all $(u_i,v_i)\in H^0$, $i=1,2$,}
\end{equation}
 and its associated norm $\|\cdot\|_{H^0}= (\cdot,\cdot)_{H^0}^{1/2}$. We also introduce the Hilbert space
 \begin{equation}\label{7}
H^1 = \left\{(u,v)\in H^1(\Omega)\times H^1(\Gamma): v=u_{|\Gamma}, v=0
\,\,\ \text{on $\Gamma_0$}\right\},
\end{equation}
endoweed with the norm inherited from the product. To simplify the notation we shall identify, when useful, $H^1$ with its isomorphic
counterpart
\begin{equation}\label{7bis}
H^1_{\Gamma_0}(\Omega,\Gamma)=\{u\in H^1(\Omega): u_{|\Gamma}\in H^1(\Gamma)\cap
L^2(\Gamma_1)\},
\end{equation}
studied for example in \cite{pucvit}, through the identification $u\mapsto (u,u_{|\Gamma})$. So we shall write, without further mention, $u\in H^1$ for
functions defined in $\Omega$.  Moreover, we shall drop the notation
$u_{|\Gamma}$, when useful, so we shall write $\int_{\Gamma_1} |u|^2$
and so on, for $u\in H^1$.

It is well--known, see \S~\ref{section 2.5} below for details, that the assumption $\mathcal{H}^{N-1}(\Gamma_0)>0$
yields a Poincar\`e type inequality in $H^1(\Omega)$. As a consequence,  the product norm on $H^1$ is equivalent to the norm $\|\cdot\|_{H^1}=(\cdot,\cdot)_{H^1}^{1/2}$ induced by the inner product
\begin{equation}\label{8}
  (u,v)_{H^1}=\int_\Omega\nabla u\nabla \overline{v}+\int_{\Gamma_1}(\nabla_\Gamma u,\nabla_\Gamma v)_\Gamma ,
\end{equation}
where $\nabla_\Gamma$ denotes the Riemannian gradient on $\Gamma$,  and $(\cdot,\cdot)_{|\Gamma}$
the unique Hermitian extension to the tangent bundle of the Riemannian metric on $\Gamma$, see \S~\ref{section 2.2}. We also denote
$|\cdot|_\Gamma =(\cdot,\cdot)|_\Gamma^{1/2}$.
\subsection{Main results} To state our main results we first make precise which type of solutions we shall consider. Given any $\lambda\in\C$, by a weak solution of \eqref{1} we shall mean $u\in H^1$ such that
\begin{equation}\label{9}
  \int_\Omega \nabla u\nabla \phi+\int_{\Gamma_1}(\nabla_\Gamma u,\nabla_\Gamma\overline{\phi})_\Gamma=
  \lambda\left(\int_\Omega u\phi+\int_{\Gamma_1}u\phi\right)\quad\text{for all $\phi\in H^1$.}
\end{equation}
Moreover, an eigenvalue for \eqref{1} will be $\lambda\in\C$ for which \eqref{1} has a nontrivial weak solution, which is called an eigenfunction. Finally, for each eigenvalue $\lambda$, the subspace of weak solutions of \eqref{1} will be called the eigenspace associated to $\lambda$. The dimension of the eigenspace is called the eigenvalue's multiplicity.

Our first result is the exact analogue  of the classical result concerning problem \eqref{4}, see \cite[\S~6.6,~Theorem~1,~p. 355]{Evans} in the real case.
\begin{thm}[\bf Spectral decomposition]\label{Theorem1}Problem \eqref{1} has countably many eigenvalues, all of which are of finite multiplicity, constituting the set $\Lambda\subset\R$.
By repeating each eigenvalue according to its multiplicity we can write $\Lambda=\{\lambda_n,n\in\N\}$, where
\begin{equation}\label{10}
0<\lambda_1\le \lambda_2\le\ldots\le\lambda_n\le \lambda_{n+1}\le\cdots,\qquad \lambda_n\to\infty\quad \text{as $n\to\infty$.}
\end{equation}
Moreover there exists an orthogonal sequence $(u_n)_n$ in $H^1$ such that each $u_n$ is a real--valued eigenfunction corresponding to $\lambda_n$, i.e. it is a nontrivial weak solution of
\begin{equation}\label{11}
\begin{cases} - \Delta u_n=\lambda_n u_n \qquad &\text{in
$\Omega$,}\\
u_n=0 &\text{on $\Gamma_0$,}\\
-\Delta_\Gamma u_n +\partial_\nu u_n =\lambda_n u_n\qquad
&\text{on
$\Gamma_1$.}
\end{cases}
\end{equation}
Moreover, $\lambda_n=\|u_n\|_{H^1}^2$ and $u_n\in C^\infty(\Omega)$.
Finally $\{u_n,n\in\N\}$ is a Hilbert basis of $H^0$ and  $\left\{\frac{u_n}{\sqrt{\lambda_n}},n\in\N\right\}$ is a Hilbert basis of $H^1$.
\end{thm}
\begin{rem}\label{Remark 1}Theorem~\ref{Theorem1} generalizes \cite[Theorem~4.4]{KnopfLiu}, in the  case
 $K=0$,  $\alpha=1$ and $\gamma=\omega$, to $\Gamma_0\not=\emptyset$ and to the complex case.
\end{rem}

The proof of Theorem~\ref{Theorem1} is based on a preliminary analysis of the doubly elliptic inhomogeneous problem
\begin{equation}\label{12}
\begin{cases} -\Delta u=f \qquad &\text{in
$\Omega$,}\\
u=0 &\text{on $\Gamma_0$,}\\
-\Delta_\Gamma u +\partial_\nu u =g\qquad
&\text{on
$\Gamma_1$,}
\end{cases}
\end{equation}
when $f\in L^2(\Omega)$ and $g\in L^2(\Gamma_1)$, see \S~\ref{section 3} below, and on
the standard Spectral Decomposition Theorem for self--adjoint compact operators.

Theorem~\ref{Theorem1} is complemented by the following variational characterization of the eigenvalues, a particular emphasis being given to the principal eigenvalue $\lambda_1$.

It extends to problem \eqref{1} well--known results concerning problem \eqref{4}, see for example \cite{Evans} or \cite{dautraylionsvol3}. 
\begin{thm}[\bf Variational characterization of the eigenvalues]\label{Theorem 2}
With the notation of Theorem~\ref{Theorem1}, we have
\begin{equation}\label{13}
\lambda_1=\min\left\{\|u\|_{H^1}^2: \quad u\in H^1,\quad
\|u\|_{H^0}=1 \right\},
\end{equation}
and the generalized Rayleigh formula
\begin{equation}\label{14}
\lambda_1=\min_{u\in H^1\setminus\{0\}}\frac {\|u\|_{H^1}^2}{\|u\|_{H^0}^2}=\min_{u\in H^1\setminus\{0\}}\frac {\int_\Omega |\nabla u|^2+\int_{\Gamma_1}|\nabla_\Gamma u|_\Gamma^2}{\int_\Omega |u|^2+\int_{\Gamma_1}|u|^2}
\end{equation}
holds true. Moreover, for any $u\in H^1$ such that $\|u\|_{H^0}=1$, $u$ is a weak solution of
\begin{equation}\label{15}
\begin{cases} - \Delta u=\lambda_1 u \qquad &\text{in
$\Omega$,}\\
u=0 &\text{on $\Gamma_0$,}\\
-\Delta_\Gamma u +\partial_\nu u =\lambda_1 u\qquad
&\text{on
$\Gamma_1$,}
\end{cases}
\end{equation}
if and only if $\|u\|_{H^1}^2=\lambda_1$.

Finally, denoting   by $\mathcal{S}_{n-1}$,  the collection
of $(n-1)$--dimensional linear subspaces of $H^1$, for $n\ge 2$, the following generalized Courant--Fischer--Weyl  min--max principle  holds:
\begin{equation}\label{16}
\lambda_n=\max_{V\in\mathcal{S}_{n-1}}\,\,\min _{u\in V^\bot, \,\|v\|_{H^0}=1}  \|u\|_{H^1}^2.
\end{equation}
\end{thm}
The generalized Rayleigh formula \eqref{14} makes evident that $1/\sqrt{\lambda_1}$ is the optimal constant of the embedding
$H^1\hookrightarrow H^0$. Moreover nontrivial weak solutions of \eqref{15} are the principal natural oscillation modes
for problem \eqref{2}.

The study of the eigenspace associated to $\lambda_1$ is then of great importance in understanding the  most important vibration behaviour of \eqref{2}. In the sequel we shall pursue this goal in a particular but important case.

Indeed we shall consider the case in which $\Gamma_0$ and $\Gamma_1$ are both relatively open on $\Gamma$, so $\Gamma$ has at least two connected components.
We remark that, in the  bass drum model described above, $\Gamma_0$ and $\Gamma_1$ are exactly these two connected components. Hence this assumption is rather natural. Moreover, we shall also assume $\Gamma$ to be at least $C^2$.

To avoid repeating these hypotheses several times we formalize them as the following assumption:
\begin{itemize}
\item[(R)]\qquad \qquad\qquad\qquad \qquad$\overline{\Gamma_0}\cap\overline{\Gamma_1}=\emptyset$ and $r\ge 2$.
\end{itemize}
When (R) holds we introduce the further spaces
\begin{alignat}2\label{17}
  &H^m=[H^m(\Omega)\times H^m(\Gamma_1)]\cap H^1, && \qquad \text{for $m\in\N$, $2\le m\le r$,} \\
  \label{18}
  &W^{2,p}= [W^{2,p}(\Omega)\times W^{2,p}(\Gamma_1)]\cap H^1, && \qquad \text{for $p\in[2,\infty)$,}
\end{alignat}
both endowed with the norms inherited from the product spaces inside square brackets.
Our third main result concerns the regularity of the eigenfunctions in Theorem~\ref{Theorem1}. Beside its independent interest, this regularity result is an essential tool in pursuing the goal described above.
\begin{thm}[\bf Regularity of the eigenfunctions]\label{Theorem3}
Let assumption (R) hold and, for all $n\in\N$, let $u_n$ be a  weak solution of \eqref{11}. Then $u_n$ enjoys  the following further regularity properties:
\renewcommand{\labelenumi}{{\roman{enumi})}}
\begin{enumerate}
\item  $u_n\in W^{2,p}$ for all $p\in [2,\infty)$, and consequently  $u_n\in C^1(\overline{\Omega})$;
\item ${u_n}\in C^2(\Gamma_1)$;
\item $u_n\in H^m$ for $m\in\N$, $2\le m\le r$, so $u_n\in C^\infty(\overline{\Omega})$ when $r=\infty$.
\end{enumerate}
\end{thm}
The proof of Theorem~\ref{Theorem3} is based on appropriate regularity results for problem \eqref{12}, see \S~\ref{section 3} below, and on  bootstrap arguments.

Our last main result brings the parallelism between  problems \eqref{1} and \eqref{4} to a probably  unexpected level. Indeed, despite of the difference between the two problems, their principal eigenspaces and eigenfunctions exhibit similar properties.
\begin{thm}[\bf The first eigenfunction]\label{Theorem4} Let assumption $(R)$ hold and let $\Omega$ be connected. Then the principal eigenvalue $\lambda_1$ in Theorem~\ref{Theorem1} is simple, i.e. $\lambda_1<\lambda_2$ in \eqref{10}, and $u_1$ has constant sign in $\Omega\cup\Gamma_1$.
\end{thm}
The most unexpected assertion in  the last statement is  that $u_1$ has constant sign also on $\Gamma_1$.
The proof Theorem~\ref{Theorem4}, indeed, uses classical arguments (see \cite{Evans}) inside $\Omega$ and \emph{ad hoc} arguments  on  $\Gamma_1$.

Finally, we would like to point out an interesting consequence of our main results in the radial case. Let
$\Omega=B_{R_2}\setminus\overline{B_{R_1}}$, $0<R_1<R_2$, and take $\Gamma_0=\partial B_{R_2}$, $\Gamma_1=\partial B_{R_1}$.
The opposite choice of $\Gamma_0$ and $\Gamma_1$ would be possible, but this one  models a circular drumhead $B_{R_2}$ with a centered circular hole $\overline{B_{R_1}}$, in the  bass drum model.
\begin{thm}[\bf The first eigenfunction in the radial case]\label{Corollary} When $\Omega$ is as described above then any weak solution of \eqref{15} is radially symmetric. Moreover, if we choose $u_1$ in Theorem~\ref{Theorem4} to be positive, then $u_1$ is strictly decreasing in the radius.
\end{thm}
The last result shows that the principal eigenfunctions behave exactly as principal eigenfunctions of problem \eqref{4} when $\Omega$ is a ball.

The paper is organized as follows. In Section~\ref{section 2} we collect the background material needed in the paper. Section~\ref{section 3} addresses the preliminary analysis of problem~\eqref{12},  which is crucial in the proofs of our main results. Finally in Section~\ref{section 4} we prove all results stated above.
\section{Preliminaries}\label{section 2}
\subsection{Notation.}\label{section 2.1}
Given a Banach space $X$, we shall denote by $I$ the identity operator, by  $X'$  its dual and by $\langle
\cdot,\cdot\rangle_X$ the duality product between them. When another Banach space $Y$ is given we shall denote by
$\mathcal{L}(X,Y)$  the space of bounded linear operators between $X$ and $Y$.

For any $p\in (1,\infty)$, we shall denote by $p'$ its H\"{o}lder conjugate, i.e. $1/p+1/p'=1$, and
for simplicity we shall denote by $\|\cdot\|_p$  the norms in $L^p(\Omega)$ and in $L^p(\Omega;\R^N)$.
Moreover, for any relatively open $\Gamma'\subset\Gamma$, we shall denote $\|\cdot\|_{p,\Gamma'}=\|\cdot\|_{L^p(\Gamma')}$.

\subsection{Riemannian operators on $\Gamma$.} \label{section 2.2}
In the sequel we shall systematically denote by $\Gamma'$, without further mention,  any relatively open subset of $\Gamma$.
  Since $\Gamma$ is of class $C^r$,  it inherits from $\R^N$ the structure of a Riemannian $C^r$ manifold, endowed with a $C^{r-1}$ Riemannian metric (see \cite{sternberg}), trivially restricting to $\Gamma'$.

In the sequel we shall use some notation of geometric nature, which is quite common when $\Gamma$ is smooth,  see \cite{Boothby, hebey, jost, taylor}. It can be easily extended to the $C^r$ case, see for example \cite{mugnvit}, and   it trivially restricts to $\Gamma'$.

We shall denote by $T(\Gamma)$ and $T^*(\Gamma)$ the tangent and cotangent bundles, standardly fiber--wise complexified (see \cite{roman}), and  by $(\cdot,\cdot)_\Gamma$ the unique Hermitian extension to $T(\Gamma)$ of the Riemannian metric inherited from $\R^N$. This Hermitian form is  given in local coordinates by $(u,v)_\Gamma=g_{ij}u^i \overline{v^j}$ for all $u,v\in T(\Gamma)$.
We notice that in the last formula we used the summation convention. We shall keep this convention in the sequel.

The metric induces the fiber--wise defined conjugate--linear Riesz isomorphisms $\flat:T(\Gamma)\to T^*(\Gamma)$ and $\sharp=\flat^{-1}:T^*(\Gamma)\to T(\Gamma)$, known as musical isomorphisms in the real case. They are defined  by the formula
$\langle \flat u,v\rangle_{T(\Gamma)}=(v,u)_\Gamma$ for $u,v\in T(\Gamma)$,
where $\langle \cdot,\cdot\rangle_{T(\Gamma)}$ denotes the fiber-wise defined duality pairing.
Hence, in local coordinates,
\begin{equation}\label{19}
  \flat u=g_{ij}\overline{u^j}dx^i, \quad\text{and}\quad \sharp \alpha=g^{ij}\overline{\alpha_j}\partial_i, \qquad \text{for all $u\in T(\Gamma)$, $\alpha\in T^*(\Gamma)$,}
\end{equation}
where $(g^{ij})=(g_{ij})^{-1}$.
The induced bundle metric on $T^*(\Gamma)$, still denoted by $(\cdot,\cdot)_\Gamma$, is then defined  by
the formula $(\alpha,\beta)_\Gamma=\langle \alpha,\sharp\beta\rangle_{T(\Gamma)}$ for all $\alpha,\beta\in T^*(\Gamma)$. Hence
\begin{equation}\label{20}
  (\alpha,\beta)_\Gamma=(\sharp \beta, \sharp\alpha)_\Gamma,\qquad \text{for all $\alpha,\beta\in T^*(\Gamma)$}.
\end{equation}
By $|\cdot|_\Gamma^2= (\cdot,\cdot)_\Gamma$  we shall denote the associated bundle norms on $T(\Gamma)$ and $T^*(\Gamma)$.

Denoting by $d_\Gamma$ the standard differential on $\Gamma$, the Riemannian gradient operator $\nabla_\Gamma$ is defined  by
setting, for $u\in C^1(\Gamma)$,
$\nabla_\Gamma u=\sharp d_\Gamma \overline{u}$,
so $\nabla_\Gamma u=g^{ij}\partial_ju\partial_i$
 in local coordinates. By \eqref{20} one trivially gets that
$(\nabla_\Gamma u,\nabla_\Gamma v)_\Gamma=(d_\Gamma u,d_\Gamma v)_\Gamma$ for all $u,v\in H^1(\Gamma)$, so in the sequel the use of vectors or forms is optional.

The Laplace--Beltrami operator $\Delta_\Gamma$ is defined as a pointwise differential operator only when $\Gamma$ is $C^2$. It can be defined in a geometrically elegant way by using the Riemannian gradient operator and the Riemannian divergence operator, dealing with tensor fields on $\Gamma$, as in \cite[\S~2.3]{mugnvit}. To avoid such a procedure we adopt here a less elegant approach, as in \cite{Dresda1}.
When $\Gamma$ is $C^2$  we thus set, for any $u\in C^2(\Gamma)$,
\begin{equation}\label{21}
\Delta_{\Gamma} u=g^{-1/2}\partial_i( g^{1/2}g^{ij}\partial_j u), \quad\text{where $g=\det (g_{ij})$,}
\end{equation}
in local coordinates. The approach in \cite{mugnvit} actually shows that \eqref{21} does not depend on the coordinate system.
\subsection{Functional spaces and operators on $\Gamma$.}\label{section 2.3}
The Riemannian gradient $\nabla_\Gamma$  is still defined, using a density argument, by $\nabla_\Gamma u=\sharp d_\Gamma \overline{u}$ for all $u\in H^1(\Gamma)$. It is well known, see for example \cite[Chapter 3]{mugnvit}, that $H^1(\Gamma)$ can be equipped with the  norm $\|\cdot\|_{H^1(\Gamma)}$, equivalent to the one introduced in \cite{grisvard}, given by
\begin{equation}\label{22}
\|u\|_{H^1(\Gamma)}^2=\|u\|_{2,\Gamma}^2+\|\nabla_\Gamma u\|_{2,\Gamma}^2,\quad\text{where}\quad \|\nabla_\Gamma u\|_{2,\Gamma}^2:=\int_\Gamma |\nabla_\Gamma u|_\Gamma^2.
\end{equation}
Moreover one can replace $\Gamma$ with $\Gamma'$ in \eqref{22}. We shall apply this remark, in particular,
when $\Gamma'=\Gamma_1$ and, when assumption (R) holds, also to $\Gamma'=\Gamma_0$.

Trivially the space $H^1$ in \eqref{7} can be equivalently defined as
$$H^1=\{(u,v)\in H^1(\Omega)\times H^1_{\Gamma_0}(\Gamma): v=u_{|\Gamma}\},$$
where
\begin{equation}\label{23}
H^1_{\Gamma_0}(\Gamma)=\{u\in H^1(\Gamma): u=0\quad\text{a.e. on $\Gamma_0$}\}.
\end{equation}
Trivially $H^1_{\Gamma_0}(\Gamma)$ is a closed subspace of $H^1(\Gamma)$. We shall endow it  with the  inherited norm.
Due to the relevance of the space $H^1$ in our analysis it is  useful to make some remarks on it.

Since, for all $u\in H^1_{\Gamma_0}(\Gamma)$, one has
$\nabla_\Gamma u=0$ a.e. on $\Gamma\setminus\overline{\Gamma_1}$, and $\mathcal{H}^{N-1}(\overline{\Gamma_0}\cap\overline{\Gamma_1})=0$, we have
\begin{equation}\label{24}
\|u\|_{H^1(\Gamma)}=\|u_{|\Gamma_1}\|_{H^1(\Gamma_1)}\qquad\text{for all $u\in H^1_{\Gamma_0}(\Gamma)$.}
\end{equation}
Formula \eqref{24} suggests the possibility of identifying $H^1_{\Gamma_0}(\Gamma)$ with $H^1(\Gamma_1)$. Actually
two different geometrical situations may occur.
\renewcommand{\labelenumi}{{\Alph{enumi})}}
\begin{enumerate}
\item One can have $\overline{\Gamma_0}\cap\overline{\Gamma_1}=\emptyset$, this case occurring when assumption (R) holds. In it, since the characteristic functions $\chi_{\Gamma_i}$, $i=0,1$ are $C^r$ on $\Gamma$, by identifying the elements of $H^1(\Gamma_i)$, $i=0,1$, with their trivial extensions to $\Gamma$, we have the splitting
    \begin{equation}\label{25}
      W^{s,p}(\Gamma)=W^{s,p}(\Gamma_0)\oplus W^{s,p}(\Gamma_1),\qquad\text{for $1<p<\infty$, $s\in\R$, $|s|\le r$.}
    \end{equation}
In accordance with the identification \eqref{5} we then have
\begin{equation}\label{26}
W^{s,p}(\Gamma_1)=\{u\in   W^{s,p}(\Gamma): u=0\quad\text{a.e. on $\Gamma_0$}\},
\end{equation}
so in particular  $H^1_{\Gamma_0}(\Gamma)=H^1(\Gamma_1)$.

\item  One can have $\overline{\Gamma_0}\cap\overline{\Gamma_1}\not=\emptyset$. In this case the set $\Gamma_0$ is not relatively open on $\Gamma$ and \eqref{25}--\eqref{26} do not hold. Indeed, for example, $\chi_{\Gamma_1}\not\in H^1_{\Gamma_0}(\Gamma)$.
In this case the elements of $H^1_{\Gamma_0}(\Gamma)$ ``vanish''  at the relative boundary $\partial\Gamma_1=\overline{\Gamma_0}\cap\overline{\Gamma_1}$ of $\Gamma_1$ on $\Gamma$, although such a notion can be made more precise only when $\partial\Gamma_1$ is regular enough. For example, if $\Gamma$ is smooth and $\overline{\Gamma_1}$ is a manifold with boundary $\partial\Gamma_1$,   see \cite[\S 5.1]{taylor}, then $H^1_{\Gamma_0}(\Gamma)$ is isometrically isomorphic to  the space
    $H^1_0(\Gamma_1):=\overline{C^\infty_c(\Gamma_1)}^{\|\cdot\|_{H^1(\Gamma_1)}}$.
\end{enumerate}
In the present paper we shall simultaneously deal with both cases A) and B) above up to the point where such a procedure is possible, that is up to the proof of Theorem~\ref{Theorem 2}. After it we shall restrict to the case A).

 Since $g$, $g^{ij}$ are $C^{r-1}$ and $\Gamma$ is compact, when $\Gamma$ is $C^2$ formula \eqref{21} extends by density to $u\in W^{2,p}(\Gamma)$, $1<p<\infty$, so defining an operator $\Delta_\Gamma\in\mathcal{L}(W^{2,p}(\Gamma);L^p(\Gamma))$. This operator restricts,  for $1<p<\infty$,  and $s\in\R$, $1\le s\le r-1$, to  $\Delta_\Gamma\in\mathcal{L}(W^{s+1,p}(\Gamma);W^{s-1,p}(\Gamma))$. Again using  the compactness of $\Gamma$  and \eqref{21}, integrating by parts and introducing a $C^2$ partition of the unity, one gets 
\begin{equation}\label{27}
-\int_\Gamma \Delta_\Gamma u \,\overline{v}=\int_\Gamma (\nabla_\Gamma u,\nabla_\Gamma v)_\Gamma\quad\text{for $u\in W^{2,p}(\Gamma)$, $v\in W^{1,p'}(\Gamma)$, $1<p<\infty$.}
\end{equation}
Formula \eqref{27} motivates the following definition of the operator 
$$\Delta_\Gamma\in\mathcal{L}(W^{1,p}(\Gamma);W^{-1,p}(\Gamma))$$
 also when $\Gamma$ is merely $C^1$. Indeed, recalling that
\begin{equation}\label{28}
W^{-s,p}(\Gamma)\simeq [W^{s,p'}(\Gamma)]',\qquad\text{for $1<p<\infty$, $s\in\R$, $0\le s\le r$,}
\end{equation}
we can set
\begin{equation}\label{29}
\langle -\Delta_\Gamma u,v\rangle _{W^{1,p'}(\Gamma)}=\int_\Gamma (\nabla_\Gamma u, \nabla_\Gamma \overline{v})_\Gamma\quad\text{for all $u\in W^{1,p}(\Gamma)$, $v\in W^{1,p'}(\Gamma)$.}
\end{equation}
By density, when $\Gamma$ is $C^2$, the operator defined in \eqref{29} is the unique  extension of
the one defined above. Hence, by interpolation, we get that, in the general case $r\ge 1$, we have
$\Delta_\Gamma\in\mathcal{L}(W^{s+1,p}(\Gamma);W^{s-1,p}(\Gamma))$ whenever $1<p<\infty$, $s\in\R$, $0\le s\le r-1$.  Using \eqref{27}--\eqref{28} the operator $\Delta_\Gamma$ can be extended, by transposition, to get
\begin{equation}\label{30}
\Delta_\Gamma\in\mathcal{L}(W^{s+1,p}(\Gamma);W^{s-1,p}(\Gamma))\quad\text{for $1<p<\infty$, $s\in\R$, $|s|\le r-1$.}
\end{equation}
 In the case A) defined above both $\Gamma_0$ and $\Gamma_1$ are relatively open and compact, so we can repeat previous arguments  and get \eqref{27} and \eqref{28} with $\Gamma$ replaced by $\Gamma_0$ and $\Gamma_1$. Hence \eqref{27} and \eqref{28} continue to hold when replacing $\Gamma$ with  $\Gamma_0$ and $\Gamma_1$, so we can set $\Delta_{\Gamma_i}\in\mathcal{L}(W^{1,p}(\Gamma_i);W^{-1,p}(\Gamma_i))$, $i=0,1$, by replacing $\Gamma$ with  $\Gamma_i$  in \eqref{29}. In this way we get
\begin{equation}\label{30bis}
\Delta_{\Gamma_i}\in\mathcal{L}(W^{s+1,p}(\Gamma_i);W^{s-1,p}(\Gamma_i))\,\,i=0,1,\, 1<p<\infty, s\in\R, |s|\le r-1.
\end{equation}
Moreover, using the splitting \eqref{25}, one obtains $\Delta_\Gamma=(\Delta_{\Gamma_0},\Delta_{\Gamma_1})$, so by the identification \eqref{26} one gets $\Delta_{\Gamma_1}={\Delta_\Gamma}_{|W^{s+1,p}(\Gamma_1)}$. Hence the symbols $\Delta_{\Gamma_1}$ and $\Delta_\Gamma$ can be used with the same meaning in the present paper.

In the sequel we shall use, only  when $p=2$ or $s=1$, the following isomorphism properties of the operator $-\Delta_{\Gamma_1}+I$.
\begin{lem}\label{Lemma 1}When $\overline{\Gamma_0}\cap\overline{\Gamma_1}=\emptyset$ the operator $$A_{s,p}=-\Delta_{\Gamma_1}+I\in\mathcal{L}(W^{s+1,p}(\Gamma_1);W^{s-1,p}(\Gamma_1))$$
 is an algebraic and topological isomorphism for $1<p<\infty$ and $s\in\R$, $|s|\le r-1$.
\end{lem}
\begin{proof} The result is well--known when $\Gamma$ is smooth, see for example \cite[p. 28]{taylor3}. Moreover the proofs in the cases $p=2$ and $s=1$ have been already explicitly given, respectively see  \cite[Theorem~5.0.1]{mugnvit} and \cite[Lemma B.1,~Appendix~B]{Dresda1}.
Since in the present paper we shall use only these cases, here we just sketch how to generalize the arguments used in the proof of
\cite[Lemma B.1,~Appendix~B]{Dresda1} from $|s|\le 1$ to $|s|\le r-1$ (this generalization being needed only when  $r>2$).

When $s\in[1,r-1]$, using the standard localization technique exactly as in the proof of \cite[Lemma~B.2,~Appendix~B]{Dresda1},
as well as elliptic higher  regularity results, see for example \cite[Chapter~2,~Theorem~2.5.1.1.,~p. 128]{grisvard}, one gets that  $-\Delta_{\Gamma_1}u+u\in W^{s-1,p}(\Gamma_1)$  implies  $u\in W^{s+1,p}(\Gamma_1)$, so $A_{s,p}$ is surjective. Being injective when $s=1$ one then obtains that $A_{s,p}$ is bijective for $s\in [0,r-1]$. The proof can then be completed by transposition.
\end{proof}

\subsection{The trace and normal derivative operators}\label{section 2.4} By  \cite[Chapter~1,~Theorem~1.5.1.2,~p.~38]{grisvard}, the standard trace operator $u\mapsto u_{|\Gamma}$ from $C(\overline{\Omega})$ to $C(\Gamma)$, when restricted to $C(\overline{\Omega})\cap W^{m,p}(\Omega)$, $1<p<\infty$, $m\in\N$, $1\le m\le r$, has a  unique surjective extension $\Tr_\Gamma\in {\mathcal L}\left(W^{m,p}(\Omega),W^{m-\frac 1p,p}(\Gamma)\right)$. Moreover, when $m=1$, the operator $\Tr_\Gamma$ has a bounded right--inverse, i.e. an operator $\mathcal{R}\in {\mathcal L}\left(W^{1-\frac 1p,p}(\Gamma),W^{1,p}(\Omega)\right)$ such that $\Tr_\Gamma\cdot \mathcal{R}=I$, with $\mathcal{R}u$ independent on $p$.
For the sake of simplicity we shall denote, as in \S~\ref{section 1}, $\Tr_\Gamma u=u_{|\Gamma}$. 

Moreover, denoting by $u_{|\Gamma_i}$, for $i=0,1$, the restriction of $u_{|\Gamma}$ to $\Gamma_i$,
and also denoting $\Tr_{\Gamma_i}u=u_{|\Gamma_i}$, one  has $\Tr_{\Gamma_1}\in {\mathcal L}\left(W^{m,p}(\Omega),W^{m-\frac 1p,p}(\Gamma_1)\right)$. When $\overline{\Gamma_0}\cap\overline{\Gamma_1}=\emptyset$ one also gets $\Tr_{\Gamma_0}\in {\mathcal L}\left(W^{m,p}(\Omega),W^{m-\frac 1p,p}(\Gamma_0)\right)$.

Moreover, when $r\ge 2$ (i.e. $\Gamma$ is $C^2$),  for all $p\in (1,\infty)$,  $m\in\N$ such that $2\le m\le r$,
and $u\in W^{m,p}(\Omega)$, denoting $\partial_\nu u =\sum_{i=1}^N \partial_i u_{|\Gamma}\,\nu^i$, where $\nu=(\nu^1,\ldots,\nu^N)$, and $\partial_i u_{|\Gamma_i}$ are taken in the trace sense sense,
we get the normal derivative operator
\begin{equation}\label{31}
 \partial_\nu\in\mathcal{L}\left(W^{m,p}(\Omega),W^{m-1-\frac 1p,p}(\Gamma)\right).
\end{equation}
Clearly $\partial_\nu u$  (its  restrictions to $\Gamma_i$ will be denoted by $\partial_\nu u_{|\Gamma_i}$) can be defined in such a trace sense only when $r\ge 2$ and $u\in W^{2,p}(\Omega)$. Hence, when $r=1$, we set
the normal derivative in a distributional sense as follows.

For any $u\in W^{1,p}(\Omega)$ such that $\Delta u\in L^p(\Omega)$ in the sense of distributions, we set, using \eqref{28}, $\partial_\nu u\in W^{-1/p,p}(\Gamma)$ by
\begin{equation}\label{32}\langle \partial_\nu u
,\psi\rangle_{W^{1-1/p',p'}(\Gamma)}=\int_\Omega
\Delta u\,\,\mathcal{R}\psi+\int_\Omega \nabla u\nabla
(\mathcal{R}\psi)
\end{equation}
for all $\psi\in W^{1-1/p',p'}(\Gamma)$, where $\mathcal{R}$ is the operator defined above. The so--defined operator
$\partial_\nu$ is linear and bounded from
$D_p(\Delta)=\{u\in W^{1,p}(\Omega): \Delta u\in
L^p(\Omega)\}$, equipped with the graph norm, to
$W^{-1/p,p}(\Gamma)$. Moreover, since for any $\Psi\in
W^{1,p'}(\Omega)$ such that $\Psi_{|\Gamma}=\psi$ we have
$\Psi-\mathcal{R}\psi\in W^{1,p'}_0(\Omega)$,  formula \eqref{32} extends
to
\begin{equation}\label{33}\langle \partial_\nu u
,\psi\rangle_{W^{1-1/p',p'}(\Gamma)}=\int_\Omega
\Delta u\Psi+\int_\Omega \nabla u\nabla\Psi\end{equation}
for all such  $\Psi$.

 Next, when $\overline{\Gamma_0}\cap\overline{\Gamma_1}=\emptyset$, by
the splitting \eqref{25}, for all $u\in W^{1,p}(\Omega)$ and $\psi\in W^{1-1/p',p'}(\Gamma)$ we have
$$\partial_\nu u=\partial_\nu u_{|\Gamma_0}+\partial_\nu u_{|\Gamma_1},\quad\text{and}\quad
\psi=\psi_{|\Gamma_0}+\psi_{|\Gamma_1},$$
where $\partial_\nu u_{|\Gamma_i}\in W^{-1/p,p}(\Gamma_i)$ and $\psi_{|\Gamma_i}\in W^{1-1/p',p'}(\Gamma_i)$ for   $i=0,1$, and
$$
\langle \partial_\nu u
,\psi\rangle_{W^{1-1/p',p'}(\Gamma)}=\sum_{i=0}^1\langle
\partial_\nu u_{|\Gamma_i}
,\psi_{|\Gamma_i}\rangle_{W^{1-1/p',p'}(\Gamma_i)}.
$$
 Hence, in particular, by \eqref{33},
\begin{equation}\label{34}
\langle \partial_\nu u_{|\Gamma_1}
,\psi\rangle_{W^{1-1/p',p'}(\Gamma_1)}=\int_\Omega
\Delta u\Psi+\int_\Omega \nabla u\nabla\Psi
\end{equation}
for all $\psi\in W^{1-1/p',p'}(\Gamma_1)$ and all $\Psi\in
W^{1,p'}(\Omega)$ such that $\Psi_{|\Gamma}=\psi$.

Clearly, when $r\ge 2$ and $u\in W^{2,p}(\Omega)$, integrating by parts using \cite[Theorem~18.1, p.~592]{LeoniSobolev2}, the so--defined normal derivative $\partial_\nu u$
 coincides with the one given by \eqref{31}, that is $\partial_\nu u\in  W^{2-1/\rho,\rho}(\Gamma)$, and   $\partial_\nu u_{|\Gamma_1}$ coincides with  its restriction to $\Gamma_1$,
 that is
$\partial_\nu u_{|\Gamma_i}\in  W^{2-1/\rho,\rho}(\Gamma_1)$.

\subsection{The space $H^1$.} \label{section 2.5} We recall, see \cite[Lemma 1, p. 2147]{vazvitHLB} (the quoted result trivially extends to $\Gamma$ of class $C^1$) that the space
$$H^1(\Omega;\Gamma)=\{(u,v)\in H^1(\Omega)\times H^1(\Gamma): v=u_{|\Gamma}\},$$
with the topology inherited from the product, can be identified with the space $\{u\in H^1(\Omega): u_{|\Gamma}\in H^1(\Gamma)\}$ and equivalently equipped with the norm $\|\cdot\|_{H^1(\Omega,\Gamma)}$ given by
$$\|u\|^2_{H^1(\Omega,\Gamma)}=\|\nabla u\|_2^2+\|\nabla_\Gamma u\|_{2,\Gamma}^2+\|u\|_{2,\Gamma}^2.$$
The identification made in \S~\ref{section 1} between the spaces $H^1$ and $H^1_{\Gamma_0}(\Omega,\Gamma)$, defined in \eqref{7} and \eqref{7bis}, is a simple consequence of the identification above. Moreover, by \eqref{24},
 $H^1$ can be equivalently equipped with the norm $\vertiii\cdot_{H^1}$ given by
$$
 \vertiii u_{H^1}^2= \|\nabla u\|_2^2+\|\nabla_\Gamma u\|_{2,\Gamma_1}^2+\|u\|_{2,\Gamma_1}^2.
$$
On the other hand, by the assumption $\mathcal{H}^{N-1}(\Gamma_0)>0$, the norm $\|\cdot\|_{H^1}$ given by \eqref{8} is equivalent to $\vertiii\cdot_{H^1}$, as it is well--known. We refer to \cite[Lemma~1,~p.~8]{MP} for an explicit proof. Hence in the paper, as stated in \S~\ref{section 1}, we shall endow the space $H^1$ with the norm $\|\cdot\|_{H^1}$ induced by the inner product $(\cdot,\cdot)_{H^1}$ defined in \eqref{8}.
\section{Analysis of problem \eqref{12}}\label{section 3}
This section is devoted to a preliminary analysis of problem \eqref{12}. We start by making precise which types of solutions we are  going to consider.
\begin{definition} Let $f\in L^2(\Omega)$ and $g\in L^2(\Gamma_1)$.
\renewcommand{\labelenumi}{{\roman{enumi})}}
\begin{enumerate}
\item We say that $u\in H^1$ is a \emph{weak solution} of \eqref{12} provided
\begin{equation}\label{35}
\int_\Omega \nabla u\nabla \phi+\int_{\Gamma_1}(\nabla_\Gamma u,\nabla_\Gamma \overline{\phi})_\Gamma=\int_\Omega f\phi+\int_{\Gamma_1}g\phi\quad\text{for all $\phi\in H^1$.}
\end{equation}
\item When assumption (R) holds, we say that $u\in H^2$ is a \emph{strong solution} of \eqref{12} provided
$$
  -\Delta u=f\quad\text{in $L^2(\Omega)$,}\quad-\Delta_{\Gamma_1}u_{|\Gamma_1}+\partial_\nu u_{|\Gamma_1}=g\quad \text{in $L^2(\Gamma_1)$,}
$$
where $\Delta_{\Gamma_1}u$ was defined in \S~\ref{section 2.3} while  $u_{|\Gamma_1}$ and  $\partial_\nu u_{|\Gamma_1}$
were defined in \S~\ref{section 2.4}.
\end{enumerate}
\end{definition}

Essentially as in \cite{Dresda1}, it is useful to deal with weak solutions of \eqref{12} in a more abstract sense.
  By \cite[Lemma~2.1]{Dresda1}, trivially extended to the complex case, the embedding $H^1\hookrightarrow H^0$ is dense.
In the sequel we shal identify $L^2(\Omega)$ and  $L^2(\Gamma_1)$ with their duals $[L^2(\Omega)]'$ and $[L^2(\Gamma_1)]'$, coherently with identifications usually made in the distribution sense. We shall also identify $L^2(\Omega)$ and  $L^2(\Gamma_1)$ with their isometric copies
$L^2(\Omega)\times\{0\}$ and  $\{0\}\times L^2(\Gamma_1)$ contained in $H^0=L^2(\Omega)\times L^2(\Gamma_1)$. Consequently, we shall identify $H^0$ with its dual $(H^0)'$, according to the identity
\begin{equation}\label{36}
  \langle u,v\rangle_{H^0}=(u,\overline{v})_{H^0}\qquad\text{for all $u,v\in H^0$.}
\end{equation}
We then introduce the chain of dense embeddings, or Gel'fand triple,
\begin{equation}\label{37}
  H^1\hookrightarrow H^0\simeq (H^0)'\hookrightarrow (H^1)',
\end{equation}
in which \eqref{36} particularizes to
\begin{equation}\label{38}
  \langle u,v\rangle_{H^1}=(u,\overline{v})_{H^0}\qquad\text{for all $u\in H^0$ and $v\in H^1$.}
\end{equation}
 Also recalling \eqref{8}, we now introduce the operator $A_1\in\mathcal{L}\left(H^1;(H^1)'\right)$ given by
\begin{equation}\label{39}
\langle A_1u,v\rangle _{H^1}=(u,v)_{H^1}=\int_\Omega\nabla u\nabla v+\int_{\Gamma_1}(\nabla_\Gamma u,\nabla_\Gamma \overline{v})_\Gamma,
\end{equation}
and its part $A:D(A)\subset H^1\subset H^0\to H^0$ given by
\begin{equation}\label{40}
D(A)=\{u\in H^1: A_1u\in H^0\},\quad Au=A_1u\quad\text{for all $u\in D(A)$},
\end{equation}
were \eqref{37} was used. By a quick comparison between \eqref{35} and \eqref{39}--\eqref{40} one gets the following result.
\begin{lem}\label{Lemma 2}Any $u\in H^1$ is a weak solution of \eqref{12} if and only if $u\in D(A)$ and $Au=(f,g)$.
\end{lem}
The following result  shows that, when assumption (R) holds, the concepts of weak and strong solution of \eqref{12} coincide.
\begin{lem}\label{Lemma 3}
If assumption (R) holds then $D(A)=H^2$ and $u$ is a strong solution of \eqref{12} if and only if it is a weak solution of it.
\end{lem}
\begin{proof}Let $u\in H^2$ be a strong solution of \eqref{12}. Multiplying the equation $-\Delta u=f$ by $\phi\in H^1$ and integrating by parts in $\Omega$  we get
\begin{equation}\label{40bis}
\int_\Omega \nabla u\nabla \phi+\int_{\Gamma_1}\partial_\nu u \phi=\int_\Omega f\phi.
\end{equation}
Moreover, multiplying the equation $-\Delta_{\Gamma_1}u_{|\Gamma_1}+\partial_\nu u_{|\Gamma_1}=g$ by $\phi$, integrating on $\Gamma_1$ and using the (integration by parts) formula \eqref{27} on $\Gamma_1$, we get
$$\int_{\Gamma_1} (\nabla_\Gamma  u,\nabla_\Gamma \overline{\phi})_\Gamma+\int_{\Gamma_1}\partial_\nu u \phi=\int_{\Gamma_1} g\phi.$$
Combining it with \eqref{40bis} we get \eqref{35}, that is  $u$ is a weak solution of \eqref{12}. By Lemma~\ref{Lemma 2} thus $u\in D(A)$. Since all $u\in H^2$ are strong solutions of \eqref{12} for appropriate $(f,g)\in H^0$, we also get $H^2\subseteq D(A)$.

Conversely, let $u\in D(A)$ be a weak solution of \eqref{12}. We first claim that $u\in H^2$, so proving that $D(A)\subseteq H^2$ and hence $D(A)= H^2$. Taking $v\in \mathcal{D}(\Omega)=C^\infty_c(\Omega)$ in \eqref{35} we first get that $-\Delta u=f$ in $\mathcal{D}'(\Omega)$. Since $f\in L^2(\Omega)$, see \S~\ref{section 2.4}, $u$ has a distributional normal derivative
$\partial_\nu u_{|\Gamma_1}\in H^{-1/2}(\Gamma_1)$ and, by \eqref{34}, we can rewrite \eqref{35} as
$$\langle \partial_\nu u_{|\Gamma_1},\phi\rangle_{H^{1/2}(\Gamma_1)}+\int_{\Gamma_1} (\nabla_\Gamma u,\nabla_\Gamma\overline{\phi})_\Gamma=\int_{\Gamma_1}g\phi\qquad\text{for all $\phi\in H^1$.}$$
By the surjectivity of the trace operator, see \S~\ref{section 2.4},  the last formula holds true for all $\phi\in H^1(\Gamma_1)$.
Hence, by \eqref{29}, we get  that $\partial_\nu u_{|\Gamma_1}-\Delta_{\Gamma_1}u=g$ in $H^{-1}(\Gamma_1)$. As a consequence of the last equation we thus get $-\Delta_{\Gamma_1}u\in H^{-1/2}(\Gamma_1)$. By Lemma~\ref{Lemma 1} we then obtain that
$u_{|\Gamma_1}\in H^{3/2}(\Gamma_1)$. Recalling that $-\Delta u=f$ in $\mathcal{D}'(\Omega)$, by elliptic regularity (see \cite[Chapter~2,~Theorem~2.4.2.5,~p. 124]{grisvard}) we then get $u\in H^2(\Omega)$. Using \eqref{31} we then have $\partial_\nu u_{|\Gamma_1}\in H^{1/2}(\Gamma_1)$, so $-\Delta_{\Gamma_1}u=g-\partial_\nu u_{|\Gamma_1}\in L^2(\Gamma_1)$.  Using Lemma~\ref{Lemma 1} again then $u_{|\Gamma_1}\in H^2(\Gamma_1)$, so $u\in H^2$ and our claim is proved.
Moreover, since $u\in H^2$, the equations $-\Delta u=f$ and $\partial_\nu u_{|\Gamma_1}-\Delta_{\Gamma_1}u=g$ respectively hold in $L^2(\Omega)$ and $L^2(\Gamma_1)$, so  we also get that $u$ is a strong solution of \eqref{12}.
\end{proof}
Before stating our next result we remark that, since the operator $A$ is trivially closed, we can endow  $D(A)$ with the graph norm $\|\cdot\|_{D(A)}=\|\cdot\|_{H^1}+\|A(\cdot)\|_{H^0}$, obtaining a Hilbert space.  The following result shows that problem \eqref{12} is well--posed and $A$ is an isomorphism.
\begin{thm}[\bf Well--posedness for problem \eqref{12}]\label{Proposition 1} For all $f\in L^2(\Omega)$ and $g\in L^2(\Gamma_1)$ problem \eqref{12} has a unique weak solution $u\in H^1$.  Moreover there is a positive constant $c_1=c_1(\Omega,\Gamma_1)$ such that
\begin{equation}\label{41}
  \|u\|_{H^1}\le c_1\left(\|f\|_2+\|g\|_{2,\Gamma_1}\right)\quad\text{for all $f\in L^2(\Omega)$ and $g\in L^2(\Gamma_1)$.}
\end{equation}
Hence $A$ is an algebraic and topological isomorphism between $D(A)$ and $H^0$, with inverse $A^{-1}\in\mathcal{L}(H^0,D(A)) \hookrightarrow\mathcal{L}(H^0,H^1)$.
\end{thm}
\begin{proof}Using \eqref{8} and \eqref{36} we can rewrite \eqref{35} as $(u,\phi)_{H^1}=((f,g),\phi)_{H^0}$ for all $\phi\in H^1$. Hence, by the complex version of the Riesz Representation Theorem (see for example
\cite[Chapter IV,~\S~6.4,~pp.~302--303]{dautraylionsvol2}), for all $(f,g)\in H^0$ problem \eqref{12} has a unique weak solution. Consequently, by Lemma~\ref{Lemma 2}, the operator $A$ is bijective from $D(A)$ onto $H^0$. Trivially $A\in\mathcal{L}(D(A);H^0)$, so by the Closed Graph Theorem we get $A^{-1}\in\mathcal{L}(H^0;D(A))\hookrightarrow\mathcal{L}(H^0,H^1)$, also proving \eqref{41}. By the way \eqref{41} is also a direct consequence of the Riesz Theorem.
\end{proof}
Our final result on problem \eqref{12} concerns  regularity properties when  assumption (R) holds.
\begin{thm}[\bf Regularity for problem \eqref{12}] \label{Theorem 6} Let assumption (R) hold and $u$ be a weak solution of \eqref{12}. \renewcommand{\labelenumi}{{\roman{enumi})}}
\begin{enumerate}
\item If $f\in L^p(\Omega)$ and $g\in L^p(\Gamma_1)$ for $p\in [2,\infty)$, then $u\in W^{2,p}$. Moreover, there is a positive constant $c_2=c_2(\Omega,\Gamma_1,p)$ such that
    \begin{equation}\label{42}
    \|u\|_{W^{2,p}}\le c_2\left(\|f\|_p+\|g\|_{p,\Gamma_1}\right)
    \end{equation}
    for all $f\in L^p(\Omega)$ and $g\in L^p(\Gamma_1)$.
\item If $f\in H^{m-2}(\Omega)$ and $g\in H^{m-2}(\Gamma_1)$ for $m\in\N$, $2\le m\le r$,  then $u\in H^m$. Moreover, there is a positive constant $c_3=c_3(\Omega,\Gamma_1,m)$ such that
    \begin{equation}\label{43}
    \|u\|_{H^m}\le c_3\left(\|f\|_{H^{m-2}(\Omega)}+\|g\|_{H^{m-2}(\Gamma_1)}\right)
    \end{equation}
    for all $f\in H^{m-2}(\Omega)$ and $g\in H^{m-2}(\Gamma_1)$.
\end{enumerate}
\end{thm}
\begin{proof} We start by proving the first assertion in ii). We argue by induction on $m$.  When  $m=2$, by Lemma~\ref{Lemma 3} we have $u\in H^2$. Now we suppose by the induction hypotesis that $m\ge 3$, the assertion holds for $m-1$, $f\in H^{m-2}(\Omega)$ and $g\in H^{m-2}(\Gamma_1)$. We claim that $u\in H^m$. By the induction hypotesis we have $u\in H^{m-1}$. Consequently, by \eqref{31}, $\partial_\nu u_{|\Gamma_1}\in H^{m-5/2}(\Gamma_1)$. Since, by Lemma~\ref{Lemma 3}, we have $\partial_\nu u_{|\Gamma_1}-\Delta_{\Gamma_1}u=g$ on $\Gamma_1$, by Lemma~\ref{Lemma 1} we get $u_{|\Gamma_1}\in H^{m-1/2}(\Gamma_1)$. By the same Lemma we also have $-\Delta u=f$ in $\Omega$, so by  elliptic higher regularity results
(see \cite[Chapter~2,~Theorem~2.5.1.1,~p. 128]{grisvard}) we get $u\in H^m(\Omega)$. A further application of \eqref{31} then yields
$\partial_\nu u_{|\Gamma_1}\in H^{m-3/2}(\Gamma_1)$. Since $\partial_\nu u_{|\Gamma_1}-\Delta_{\Gamma_1}u=g$ on $\Gamma_1$, again by Lemma~\ref{Lemma 1} we get $u_{|\Gamma_1}\in H^m(\Gamma_1)$, proving our claim.

To prove  \eqref{43} we now remark that, when assumption (R) holds, the operator $A$ defined in \eqref{40} can be  rewritten in a more explicit form. Indeed, recalling that $H^0=L^2(\Omega)\times L^2(\Gamma_1)$, we can rewrite it  as the operator $A\in \mathcal{L}(H^2,H^0)$ given by $Au=(-\Delta u, -\Delta_{\Gamma_1}u_{|\Gamma_1}+\partial_\nu u_{|\Gamma_1})$. Since, trivially,
$A\in \mathcal{L}(H^m,H^{m-2})$ for all $m\in\N$, $2\le m\le r$, the first assertion shows that this restriction of $A$ is bijective, hence \eqref{43} follows by the Closed Graph Theorem.

We now turn to proving i), starting from the first assertion. When $p=2$ there is nothing to prove, so we take
$p\in (2,\infty)$, $f\in L^p(\Omega)$, $g\in L^p(\Gamma_1)$ and we claim that $u\in W^{2,p}$. Since, by Lemma~\ref{Lemma 3}, we have $u\in H^2(\Omega)$, using \eqref{31} we get $\partial_\nu u_{|\Gamma_1}\in H^{1/2}(\Gamma_1)$.
We remark that $H^{1/2}(\Gamma_1)$ coincides with the Besov space  $B^{1/2,2}_2(\Gamma_1)$, as proved in \cite[pp.~189--190]{triebel}.
We now distinguish between the cases $N=2$ and $N\ge 3$.

When $N=2$, we apply the Sobolev Embedding Theorem for Besov spaces in the critical case,
 see  \cite[Chapter~17, Theorem 17.55, p. 564]{LeoniSobolev2}, and we get $\partial_\nu u_{|\Gamma_1}\in L^q(\Gamma_1)$ for all $q\in [2,\infty)$, so $\partial_\nu u_{|\Gamma_1}\in L^p(\Gamma_1)$. Since $\partial_\nu u_{|\Gamma_1}-\Delta_{\Gamma_1}u=g$ on $\Gamma_1$, by Lemma~\ref{Lemma 1} we then obtain $u_{|\Gamma_1}\in W^{2,p}(\Gamma_1)$.
 Using  elliptic regularity again we then get $u\in W^{2,p}(\Omega)$, so $u\in W^{2,p}$, proving our claim.

 Let us now consider the case $N\ge 3$. In this case we apply the Sobolev Embedding Theorem for Besov spaces in the subcritical case, since $2<\frac{N-1}{1/2}=2(N-1)$, see  \cite[Chapter~17, Theorem 17.49, p. 561]{LeoniSobolev2}, to get that $\partial_\nu u_{|\Gamma_1}\in L^{p_1}(\Gamma_1)$, where $p_1:=\frac {2(N-1)}{N-2}>1$. If $p\le p_1$ we have $\partial_\nu u_{|\Gamma_1}\in L^p(\Gamma_1)$ and we can complete the proof of our claim as in the case $N=2$. Hence in the sequel we suppose that $p_1<p$. Using Lemma~\ref{Lemma 1} once again we then get $u_{|\Gamma_1}\in W^{2,p_1}(\Gamma_1)$ and then, by elliptic regularity, as before, $u\in W^{2,p_1}(\Omega)$. By \eqref{31} we then obtain $\partial_\nu u_{|\Gamma_1}\in W^{1-\frac 1{p_1},p_1}(\Gamma_1)$.
 We now have to distinguish between two further cases: either $p_1\ge \frac{N-1}{1-1/p_1}$ or $p_1< \frac{N-1}{1-1/p_1}$.
 The first case, which is best rewritten as $p_1\ge N$, occurs when $N=3$, while the  second one when $N\ge 4$.

 In the first case, again by the Sobolev Embedding Theorem for Besov spaces in the critical case, we get $\partial_\nu u_{|\Gamma_1}\in L^p(\Gamma_1)$ and we conclude the proof of our claim as in the case $N=2$. Hence we can suppose $N\ge 4$ in the sequel. In this case, by the Sobolev Embedding Theorem for Besov spaces in the subcritical case, we get $\partial_\nu u_{|\Gamma_1}\in L^{p_2}(\Gamma_1)$, where
 $p_2:=\frac {(N-1)p_1}{N-p_1}>p_1$. If $p\le p_2$ we get $\partial_\nu u_{|\Gamma_1}\in L^p(\Gamma_1)$ and we conclude the proof of our claim  as in the case $N=2$. Hence we take $p_2<p$ in the sequel. A further application of Lemma~\ref{Lemma 1} yields $u_{|\Gamma_1}\in W^{2,p_2}(\Gamma_1)$ and then, by elliptic regularity, as before, $u\in W^{2,p_2}(\Omega)$. By \eqref{31} then
 $\partial_\nu u_{|\Gamma_1}\in W^{1-\frac 1{p_2},p_2}(\Gamma_1)$.

It is then clear that, going on in this way, the following alternative occurs: either after finitely many iterations we get
$\partial_\nu u_{|\Gamma_1}\in L^p(\Gamma_1)$, so the proof of our claim  can be completed as when $N=2$, or we can go on indefinitely. In this case
the recursive formula $p_0=2$, $p_{n+1}=\frac {(N-1)p_n}{N-p_n}$, defines a strictly increasing sequence $(p_n)_n$ such that $1<p_n<\min\{p,N\}$ for all $n\in\N$. Let $l:=\lim_n p_n\in (1,\min\{p,N\}]$. Passing to the limit in the recursive formula  we get that $l=N$ leads to $l=\infty$, a contradiction, so $l<N$. Passing to the limit then
$l=\frac {(N-1)l}{N-l}$, that is $l=1$, a contradiction. Hence after finitely many iterations we get
$\partial_\nu u_{|\Gamma_1}\in L^p(\Gamma_1)$ and conclude the proof of our claim.
To prove \eqref{42} and conclude the proof one then applies the same argument already used to prove \eqref{43}.
\end{proof}
\section{Proofs}\label{section 4} This section is devoted to prove our main results, stated in Section~\ref{section 1}.
\begin{proof}[\bf Proof of Theorem~\ref{Theorem1}] By virtue of the Gel'fand triple \eqref{37} we can regard the operator $A^{-1}$ in Theorem~\ref{Proposition 1} as an operator $B\in\mathcal{L}(H^0)$. Moreover, since the embeddings
$H^1(\Omega)\hookrightarrow L^2(\Omega)$ and $H^1(\Gamma)\hookrightarrow L^2(\Gamma)$ are both compacts, the embedding
$H^1\hookrightarrow H^0$ is compact as well. Hence $B$ is compact. Moreover, by \eqref{36} and \eqref{38}--\eqref{39}, for all $u,v\in D(A)$ we have
$$(Au,v)_{H^0}=(A_1 u,v)_{H^0}=\langle A_1u,\overline{v}\rangle_{H^1}=(u,v)_{H^1}.$$
Hence $A$ is symmetric. Since $A$ is bijective its bounded inverse $B$ is self--adjoint.

We are then going to apply to $B$  the standard Spectral Decomposition Theorem for compact and self adjoint operators, see for example \cite{brezis2} or \cite{dautraylionsvol3}. The spectrum $\sigma(B)$ of $B$ then consists of the union of $\{0\}$ and of its point spectrum $\sigma_p(B)\subset\R$, which is at most countable. Moreover, since $A$ is bijective, $0\not\in \sigma_p(B)$ and, by the Fredholm Alternative, $\text{ker} (B-\mu I)$ is finite--dimensional for all $\mu\in\sigma_p(B)$. Moreover $H^0$, being separable and infinite dimensional, admits a Hilbert basis of eigenvectors of $B$. Hence $\sigma_p(B)=\{\mu_n, n\in\N\}$, where the sequence $(\mu_n)_n$ is injective and $\mu_n\to 0$ as $n\to\infty$. Consequently $\Lambda=\{\lambda_n,n\in\N\}$, where $\lambda_n=1/\mu_n$. Let us denote $H_n=\text{ker} (B-\mu_n I)$. Since, for each $u\in H_n$, we have $\lambda_n=\|u\|_{H^1}^2$, it follows that $\Lambda\subset (0,\infty)$. Consequently we can re--arrange the $\lambda_n$'s, according to their finite multiplicities, as in \eqref{10}.

Since, for each $n\in\N$, $H_n$ is exactly the space of weak solutions of \eqref{11}, it is invariant with respect to conjugation. Consequently, see \cite[Chapter~7,~Proposition~7.4.1]{BergerMarcel}, denoting by $H_n^\R$ the real subspace of the real--valued weak solutions of \eqref{11}, we have $H_n=H_n^\R+iH_n^\R$. Since the restriction to $H_n^\R$ of the inner product $(\cdot,\cdot)_{H^0}$ is a real scalar product, for each $n\in\N$ one can build an orthonormal real basis of it, which is also an orthonormal basis of $H_n$.
Since $H^0=\oplus_{n=1}^\infty H_n$ and the $H_n$'s are mutually orthogonal, we can then construct a Hilbert basis $\{(u_n,v_n), n\in\N\}$ of $H^0$ such that $(u_n,v_n)\in H^1$ for all $n\in\N$, so $v_n={u_n}_{|\Gamma}$ and, using the embedding $H^1\hookrightarrow H^0$ together with the identification between the spaces in \eqref{7} and in \eqref{7bis}, we shall simply write $(u_n,v_n)=u_n$.

Clearly each $u_n$ is a real--valued weak solution of \eqref{11}. By standard elliptic regularity, see \cite[Chapter~6,~Theorem~3,~p.~334]{Evans}, we have $u_n\in C^\infty(\Omega)$.
Moreover, by \eqref{8} and \eqref{9}, for all $n,m\in\N$ we have $(u_n,u_m)_{H^1}=\lambda_n (u_n,u_m)_{H^0}$, so $\lambda_n=\|u_n\|_{H^1}^2$ for all $n\in\N$ and $\left\{u_n/\sqrt{\lambda_n}, n\in\N\right\}$ is an orthonormal system on $H^1$. To recognize that it is a Hilbert basis of $H^1$, and complete the proof, we remark that, using \eqref{8} and \eqref{9} again, we have
\begin{equation}\label{caxa}
  (u,u_n)_{H^1}=\lambda_n (u,u_n)_{H^0}\qquad\text{for all $u\in H^1$.}
\end{equation}
  Hence $(u,u_n)_{H^1}=0$ for all $n\in\N$ yields $(u,u_n)_{H^0}=0$ for all $n\in\N$, and consequently  $\text{span}\{u_n,n\in\N\}$ is dense in $H^1$.
\end{proof}
\begin{proof}[\bf Proof of Theorem~\ref{Theorem 2}] We keep the notation of the previous proof. We first prove \eqref{13}--\eqref{14}.
As $\left\{u_n/\sqrt{\lambda_n}, n\in\N\right\}$ is a Hilbert basis of $H^1$, also using \eqref{caxa} we get that, for all $u\in H^1$,
\begin{equation}\label{43bis}
  u = \sum_{n=1}^\infty (u,u_n)_{H^1}\,\frac {u_n}{\lambda_n} = \sum_{n=1}^\infty (u,u_n)_{H^0} \,u_n\qquad\text{in $H^1$.}
\end{equation}
Consequently, for all $u\in H^1$ such that  $\|u\|_{H^0}=1$, we have
\begin{equation}\label{44}
 \| u\|_{H^1}^2 = \sum_{n=1}^\infty |(u,u_n)_{H^0}|^2\,\, \|u_n\|_{H^1}^2 = \sum_{n=1}^\infty \lambda_n |(u,u_n)_{H^0}|^2.
 \end{equation}
 Hence, by \eqref{10},
 \begin{equation}\label{45}
 \| u\|_{H^1}^2 \ge \lambda_1 \sum_{n=1}^\infty  |(u,u_n)_{H^0}|^2 = \lambda_1 \|u\|_{H^0}^2=\lambda_1.
 \end{equation}
 Since $\|u_1\|_{H^1}^2=\lambda_1$, the proof of \eqref{13} is complete. By a standard homogeneity argument \eqref{14}
 follows from \eqref{13}.

 To prove the second assertion in the statement let us take $u\in H^1$ with $\|u\|_{H^0}=1$. If $u$ is a weak solution of \eqref{15}, by \eqref{8} and \eqref{9} one immediately gets $\|u\|_{H^1}^2=\lambda_1$. Conversely, if $\|u\|_{H^1}^2=\lambda_1$, by \eqref{44} we get
$$\lambda_1 \sum_{n=1}^\infty  |(u,u_n)_{H^0}|^2=\lambda_1 \|u\|_{H^0}^2=\lambda_1=\|u\|_{H^1}^2=\sum_{n=1}^\infty  \lambda_n |(u,u_n)_{H^0}|^2$$
so $\sum\limits_{n=1}^\infty  (\lambda_n-\lambda_1) |(u,u_n)_{H^0}|^2=0$. Consequently, for each $n\in\N$, either $(u,u_n)_{H^0}=0$ or $\lambda_n=\lambda_1$. Denoting by $\nu_1$ the finite multiplicity of $\lambda_1$, by \eqref{43bis} we then get
$u = \sum_{n=1}^{\nu_1} (u,u_n)_{H^0} u_n$, so $u$ is a weak solution of \eqref{15}, proving our second assertion.

To complete the proof we now have to prove the generalized Courant--Fischer--Weyl formula \eqref{16}. To achieve this goal we are going to apply the min--max principle \cite[Chapter~VIII,~Theorem~10,~p.~102]{dautraylionsvol3}. We consequently have to verify the assumptions of the quoted result, that is the structural assumptions (2.73)--(2.74),  p. 39, and  (2.399), p. 98, of the quoted reference. We chose $V=H^1$ and $a=(\cdot,\cdot)_{H^1}$, so that (2.73) trivially holds. Moreover, also choosing $H=H^0$, assumption
 (2.399) is nothing but the Gel--fand triple \eqref{37}, combined with the already remarked compactness of the embedding $H^1\hookrightarrow H^0$. The quoted result then concerns the unbounded operator $\hat{A}$ defined by (2.74), that is by
 \begin{equation}\label{47}
D(\hat{A})=\{u\in H^1: v\mapsto (u,v)_{H^1}\,\text{is continuous on $H^1$ \!\!for the topology of $H^0$}\}
\end{equation}
and
\begin{equation}\label{48}
(u,v)_{H^1}=(\hat{A}u,v)_{H^0}\quad\text{for all $u\in D(\hat{A})$ and $v\in H^1$.}
\end{equation}
By comparing \eqref{47} with \eqref{40}, by virtue of the identifications in \eqref{37}, one gets $D(\hat{A})=D(A)$. Hence, by comparing \eqref{48} with \eqref{39}, one also gets $\hat{A}=A$. We can then apply the quoted result to $A$. Since \eqref{16} is nothing but \cite[Chapter~VIII,~(2.423),~p.~102]{dautraylionsvol3}, the proof is completed.
\end{proof}
\begin{proof}[\bf Proof of Theorem~\ref{Theorem3}] At first we prove i), by combining Theorem~\ref{Theorem 6}--i) with a bootstrap argument. More in detail, by Lemma~\ref{Lemma 3},  $u_n\in H^2$ for each $n\in\N$, that is  $u_n\in H^2(\Omega)$ and ${u_n}_{|\Gamma_1}\in H^2(\Gamma_1)$.
We now distinguish between the cases $N\le 4$ and $N\ge 5$. In the first one, by applying Morrey's Theorem when $N<4$ and the limiting case of Sobolev Embedding Theorem when $N=4$, we get $u_n\in L^p(\Omega)$ for all $p\in [2,\infty)$. By applying the same results on $\Gamma_1$ we get that ${u_n}_{|\Gamma_1}\in L^p(\Gamma_1)$ for all $p\in [2,\infty)$ provided $N\le 5$, so a -- fortiori when $N\le 4$. By applying Theorem~\ref{Theorem 6}--i) we then get $u_n\in W^{2,p}$ for all $p\in [2,\infty)$.
When $N\ge 5$ we further distinguish between the case $5\le N\le 8$ and $N\ge 9$. In the first one, by Sobolev Embedding Theorem, we have $u_n\in L^{q_1}(\Omega)$, where $\frac 1{q_1}=\frac 12-\frac 2N$, and as above one also gets that $u_n\in L^{q_1}(\Gamma_1)$. Hence, by Theorem~\ref{Theorem 6}--i), one has $u_n\in W^{2,q_1}$. Since $\frac 1{q_1}-\frac 2N=\frac 12-\frac 4N\le 0$, by applying
Morrey's Theorem when $N<8$ and the limiting case of Sobolev Embedding Theorem when $N=8$, we get $u_n\in L^p(\Omega)$ and ${u_n}_{|\Gamma_1}\in L^p(\Gamma_1)$ for all $p\in [2,\infty)$, so $u_n\in W^{2,p}$ for all $p\in [2,\infty)$ as in the previous case.
It is then clear how to handle the remaining case $N\ge 9$, proving after finitely many iterations that $u_n\in L^p(\Omega)$ and ${u_n}_{|\Gamma_1}\in L^p(\Gamma_1)$ , so $u_n\in W^{2,p}$, for all $p\in [2,\infty)$.
The conclusion $u_n\in C^1(\overline{\Omega})$ then follows from Morrey's Theorem.

To prove ii) we remark that, since $u_n\in W^{2,p}$, for all $p\in [2,\infty)$, using \eqref{31} we have $\partial_\nu {u_n}_{|\Gamma_1}\in W^{1-\frac 1{p},p}(\Gamma_1)$ for all $p\in [2,\infty)$.
We also remark that $W^{1-\frac 1{p},p}(\Gamma_1)$ coincides with the Besov space  $B^{1-\frac 1{p},p}_p(\Gamma_1)$, as proved in \cite[pp.~189--190]{triebel}.
By Morrey's Theorem for Besov spaces, see  \cite[Chapter~17, Theorem 17.52, p. 562]{LeoniSobolev2}, taking $p>N$, we have $\partial_\nu {u_n}_{|\Gamma_1}\in C^{0, 1-\frac Np}(\Gamma_1)$. Moreover, since we also have ${u_n}_{|\Gamma_1}\in C^1(\Gamma_1)$,  by \eqref{11} and Lemma~\ref{Lemma 3}, $-\Delta_{\Gamma_1}u_n=f_n\in C^{0, 1-\frac Np}(\Gamma_1)$. Being $\Gamma_1$ compact, we can use the
standard localization technique together with the classical Schauder Estimates, see \cite[Theorem~9.33]{brezis2} or \cite[Chapter~6]{gilbarg}, to recognize that in local coordinates $u_n$ is of class $C^{2,1-\frac Np}$, so $u_n\in C^2(\Gamma_1)$, proving ii).

To prove iii) we now take $m\in\N$, $2\le m\le r$. By a reiterate application of Theorem~\ref{Theorem 6}--ii), i.e. a bootstrap argument, one gets $u_n\in H^m$, so, by Morrey's Theorem,  $u_n\in C^\infty(\overline{\Omega})$ when $r=\infty$.
\end{proof}
\begin{proof}[\bf Proof of Theorem~\ref{Theorem4}] The proof consists in a nontrivial adaptation of the classical arguments in \cite[Chapter~6,~Proof~of~Theorem~2,~p.~356]{Evans}. At first we claim that, if $u\in H^1$ is a nontrivial real--valued weak solution of \eqref{15}, then either $u>0$ or $u<0$ in $\Omega$.

We preliminarily remark that, by Theorems~\ref{Theorem1} and \ref{Theorem3},  $u\in C^\infty(\Omega)\cap C^1(\overline{\Omega})$.
We assume, without restriction, that $\|u\|_{H^0}=1$. Denoting by $u^+$ and $u^-$ the positive and negative parts of $u$, by \eqref{6} we have
\begin{equation}\label{49}
1=\|u\|_{H^0}^2=\|u^+\|_2^2+\|u^-\|_2^2+\|u^+\|_{2,\Gamma_1}^2+\|u^-\|_{2,\Gamma_1}^2=\|u^+\|_{H^0}^2+\|u^-\|_{H^0}^2.
\end{equation}
As proved in the quoted reference, one has
\begin{equation}\label{50}
\nabla u^+=
\begin{cases}
\nabla u \quad&\text{a.e. in $\Omega^+$,}\\
    \,\,\,0 \quad&\text{a.e. in $\Omega^-$,}
\end{cases}
\qquad
\nabla u^-=
\begin{cases}
- \nabla u \quad&\text{a.e. in $\Omega^-$,}\\
    \quad \,\,0 \quad&\text{a.e. in $\Omega^+$,}
\end{cases}
\end{equation}
where $\Omega^+=\{x\in\Omega : u(x)\ge 0\}$ and $\Omega^-=\{x\in\Omega : u(x)\le 0\}$.
Since $\nabla_\Gamma u=g^{ij}\partial_ju\partial_i$  in local coordinates, by  \eqref{50} we also get
\begin{equation}\label{51}
\nabla_\Gamma u^+=
\begin{cases}
\nabla_\Gamma u \quad&\text{a.e. in $\Gamma_1^+$,}\\
    \,\,\,0 \quad&\text{a.e. in $\Gamma_1^-$,}
\end{cases}
\qquad
\nabla u^-=
\begin{cases}
- \nabla_\Gamma u \quad&\text{a.e. in $\Gamma_1^-$,}\\
    \quad \,\,0 \quad&\text{a.e. in $\Gamma_1^+$,}
\end{cases}
\end{equation}
where $\Gamma_1^+=\{x\in\Gamma_1 : u(x)\ge 0\}$ and $\Gamma_1^-=\{x\in\Gamma_1 : u(x)\le 0\}$.
Then $u^+,u^-\in H^1$ and, using \eqref{8}, \eqref{50} and \eqref{51}, we get
\begin{multline}\label{52}
 \|u\|_{H^1}^2=\int_\Omega |\nabla u|^2+\int_{\Gamma_1} |\nabla_\Gamma u|_\Gamma^2
  = \int_{\Omega^+} |\nabla u^+|^2+ \int_{\Omega^-} |\nabla u^-|^2\\ +\int_{\Gamma_1^+} |\nabla_\Gamma u^+|_\Gamma^2
+\int_{\Gamma_1^-} |\nabla_\Gamma u^-|_\Gamma^2 =\|u^+\|_{H^1}^2+\|u^-\|_{H^1}^2.
\end{multline}
By Theorem~\ref{Theorem 2} we have
\begin{equation}\label{53}
\|u^+\|_{H^1}^2\ge \lambda_1 \|u^+\|_{H^0}^2,\qquad\text{and}\quad
\|u^-\|_{H^1}^2\ge \lambda_1 \|u^-\|_{H^0}^2.
\end{equation}
Consequently, by \eqref{49}, \eqref{52} and \eqref{53},
$$ \lambda_1=\lambda_1 \left (\|u^+\|_{H^0}^2+\|u^-\|_{H^0}^2\right)\le \|u^+\|_{H^1}^2+\|u^-\|_{H^1}^2=\|u\|_{H^1}^2. $$
Since, by Theorem~\ref{Theorem 2}, we also have $\|u\|_{H^1}^2=\lambda_1$, we arrive to
$$ \lambda_1 \left (\|u^+\|_{H^0}^2+\|u^-\|_{H^0}^2\right)=\|u^+\|_{H^1}^2+\|u^-\|_{H^1}^2.$$
Combining it with \eqref{53} we obtain
\begin{equation}\label{54}
\|u^+\|_{H^1}^2 = \lambda_1 \|u^+\|_{H^0}^2,\qquad\text{and}\quad
\|u^-\|_{H^1}^2 = \lambda_1 \|u^-\|_{H^0}^2.
\end{equation}
By combining \eqref{54} and Theorem~\ref{Theorem 2} we thus get that $u^+$ and $u^-$ are weak solutions of \eqref{15}. Hence, by applying again Theorem~\ref{Theorem3}, we have $u^+,u^-\in C^\infty(\Omega)\cap C^1(\overline{\Omega})$. Since $\Omega$ is connected and $\lambda_1 u^+,\lambda_1 u^-\ge 0$ in $\Omega$, by applying the Strong Maximum Principle, see \cite[Chapter~6,~Theorem~4,~p.~350]{Evans}, we get that either $u^+\equiv 0$ or $u^+>0$ in $\Omega$, the same alternative applying to $u^-$. Since $u$ is nontrivial our claim follows.

We now claim that $\lambda_1$ is simple. By adapting the argument of the quoted reference to the complex case, let $u,w\in H^1$ be
two nontrivial weak solutions of \eqref{15}. Hence also $\text{Re}\,w$ and $\text{Im}\,w$ are real--valued weak solutions of \eqref{15}, and at least one of them is nontrivial. Hence, by our previous claim, either $\int_\Omega  \text{Re}\, w\not=0$ or $\int_\Omega  \text{Im}\, w\not=0$, so in any case $\int_\Omega w\not=0$. Consequently the equation $\int_\Omega u-\chi w=0$ exactly has a solution $\chi\in\C$.
Now, also $\text{Re}(u-\chi w)$ and $\text{Im}(u-\chi w)$ are real--valued weak solutions of \eqref{15}. Moreover, by the choice of $\chi$,  we have
$$\int_\Omega \text{Re}(u-\chi w)=\int_\Omega \text{Im}(u-\chi w)=0.$$
Our previous claim then yields $\text{Re}(u-\chi w)=\text{Im}(u-\chi w)=0$, that is $u=\chi w$ in $\Omega$, proving the current claim.

To complete the proof we now claim that, if $u\in H^1$ is a weak solution  of \eqref{15} and $u>0$ in $\Omega$, then we have $u>0$ on $\Gamma_1$ as well. Since $u\in C^1(\overline{\Omega})$ we clearly have $u\ge 0$ in $\overline{\Omega}$. Now we suppose by contradiction that there is $x_0\in \Gamma_1$ such that $u(x_0)=0$. Hence $x_0$ is an absolute minimum point for $u$ on $\Gamma_1$, from which we get $\nabla_\Gamma u(x_0)=0$.

We now recall the construction of a normal coordinate system in a neighborhood of $x_0$, made in \cite[Chapter~1,~Theorems~1.4.3~and~1.4.4]{jost}, in which $(g_{ij}(x_0))=I$. Although, in the quoted reference, $\Gamma_1$ is supposed to be a smooth manifold, we remark that the $C^2$ regularity of $\Gamma_1$ assumed here allows to repeat the construction.
Now, since $u\in C^2(\Gamma_1)$, in this coordinate system we have
\begin{equation}\label{54bis}
\Delta_\Gamma u(x_0)=g^{-1/2}\partial_i (g^{1/2}g^{ij}\partial_j u)(x_0)=  \partial^2_{ii} u(x_0)\ge 0
\end{equation}
since $x_0$ is a minimum.

We now remark that, since $u\in C^2(\Gamma_1)$ and $\partial_\nu u\in C^1(\Gamma_1)$, by Lemma~\ref{Lemma 3} the equation $-\Delta_\Gamma u+\partial_\nu u=\lambda_1 u$ on $L^2(\Gamma_1)$ actually holds in the space $C(\Gamma_1)$, i.e. pointwise. Hence, by \eqref{54bis}, we have
\begin{equation}\label{54ter}
\partial_\nu u(x_0)=\Delta_\Gamma u(x_0)+\lambda_1 u(x_0)=\Delta_\Gamma u(x_0)\ge 0.
\end{equation}

Moreover, since $\Gamma$ is $C^2$, it is also $C^{1,1}$, so $\Omega$ satisfies the interior ball condition, that is to say  there is an open ball $B\subseteq \Omega$ such that $x_0\in\partial B$. See \cite[Corollary~2,p.~550]{LewickaPeres} or \cite[Theorem~1.0.9,p.~7]{barb}. Now, since $u\in C^\infty(B)\cap C^1(\overline{B})$,  $-\Delta u=\lambda_1 u\ge 0$ in $B$, and $u(x)>u(x_0)=0$ for all $x\in B$, we can apply the classical H\"{o}pf Lemma, see \cite[Chapter~6,p.~330]{Evans}, to conclude that
$\partial_\nu u(x_0)<0$,  contradicting \eqref{54ter} and completing the proof.
\end{proof}
\begin{proof}[\bf Proof of Theorem~\ref{Corollary}]
Let $\Omega=B_{R_2}\setminus\overline{B_{R_1}}$, $0<R_1<R_2$, and $\Gamma_0=\partial B_{R_2}$, $\Gamma_1=\partial B_{R_1}$.
We first claim that, given any rotation $H:\R^N\to\R^N$, we have
\begin{equation}\label{57}
\|u\cdot H\|_{H^1}=\|u\|_{H^1} \qquad\text{for all $u\in H^1$.}
\end{equation}
Trivially, by changing variables in the multiple integral, we have
\begin{equation}\label{55}
 \|\nabla (u\cdot H)\|_2^2=\int_\Omega |\nabla (u( H(x))|^2\,dx = \int_\Omega |\nabla (u(y)|^2|\,dy=\|\nabla u\|_2^2.
\end{equation}
Moreover, since $H$ is a rotation, the metric $(\cdot,\cdot)_\Gamma$ has the same matrix representation in the spherical coordinates $(y^1,\ldots,y^{N-1})$ on $\Gamma_1=\partial B_{R_1}$, given by $\psi:\partial B_{R_1}\to \R^{N-1}$, and in the spherical coordinates on $\Gamma_1$ given by the map $\psi\cdot H:\partial B_{R_1}\to \R^{N-1}$. Hence, in these two local coordinates sytems, the volume element $d\mathcal{H}^{N-1}=\sqrt{g}dy^1\wedge\ldots\wedge dy^{N-1}$, the Riemannian gradient $\nabla_\Gamma$ and $|\cdot|_\Gamma|^2$ are equals. Hence, since $H:\Gamma_1\to\Gamma_1$ is bijective and $C^\infty$, we have
\begin{equation}\label{56}
 \|\nabla_\Gamma (u\cdot H)\|_{2,\Gamma_1}^2=\int_{\Gamma_1} |\nabla_\Gamma (u\cdot H)|_\Gamma^2 = \int_{\Gamma_1}
 |\nabla_\Gamma u|_\Gamma^2|=\|\nabla_\Gamma u\|_{2,\Gamma_1}^2.
\end{equation}
Trivially \eqref{57} follows form \eqref{55} and \eqref{56}, proving our claim.
By the same arguments one also gets that
\begin{equation}\label{58}
\|u\cdot H\|_{H^0}=\|u\|_{H^0} \qquad\text{for all $u\in H^0$.}
\end{equation}

Now let $u$ be a real--valued weak solution of \eqref{15}.  We claim that $u$ is radially symmetric. If $u$ is trivial there is nothing to prove. Hence, using Theorem~\ref{Theorem4}, we can suppose that $u$ is positive in $\Omega\cup\Gamma_1$.
By \eqref{57}, \eqref{58}, the minimality asserted in Theorem~\ref{Theorem 2} and the simplicity of the $\lambda_1$ asserted in Theorem~\ref{Theorem4}, we get that for any rotation $H$ there is $\mu=\mu_H\in \R$ such that $u\cdot H=\mu_H u$ in $\Omega$. By \eqref{58} we also have $\mu_H\in\{-1,1\}$. Since $u$ is positive then $\mu_H=1$, so $u\cdot H=u$ in $\Omega$. Since $H$ is arbitrary
our claim is proved. Passing to real and imaginary parts all weak solutions of \eqref{15} are radially symmetric.

Now let $u$ be a  positive solution of  \eqref{15}. By the previous claim we have $u(x)=w(\rho)$, $\rho=|x|$, and by Theorem~\ref{Theorem3},  $w\in C^\infty([R_1,R_2])$. Problem \eqref{15} can then be rewritten as
\begin{equation}\label{59}
\begin{cases} -w''(\rho)-\frac {N-1}\rho w'(\rho)=\lambda_1 w(\rho) \quad &\text{in $(R_1,R_2)$,}\\
-w'(R_1)=\lambda_1 w(R_1),\quad w(R_2)=0.&
\end{cases}
\end{equation}
Using a standard argument, by multiplying \eqref{59}$_1$ by $\rho^{N-1}$ one gets
$$\frac d{d\rho}[\rho^{N-1} w'(\rho)] =-\lambda_1 \rho^{N-1} w(\rho).$$
Integrating from $R_1$ to $\rho\in [R_1,R_2])$ and using the boundary conditions we consequently get
\begin{equation}\label{60}
\rho^{N-1}w'(\rho)=-\lambda_1 \int_{R_1}^\rho t^{N-1}w(t)\,dt-\lambda_1 R_1^{N-1} w(R_1).
\end{equation}
Since $w>0$ in $[R_1,R_2)$, \eqref{60} yields $w'<0$ in  $[R_1,R_2]$, completing the proof.
 \end{proof}

\bigskip

\centerline{\bf Statements and Declarations.}

\bigskip

{\bf Conflict of Interest.} The author declares that he has  no
conflict of interest.

{\bf Data availability.} Data sharing not applicable to this article as no datasets were generated or analysed during the current study.

\begin{thebibliography}{10}

\bibitem{adams}
R.~A. Adams, \emph{Sobolev spaces}, Academic Press, New York-London, 1975, Pure
  and Applied Mathematics, Vol. 65.

\bibitem{andrewskuttlershillor}
K.~T. Andrews, K.~L. Kuttler, and M.~Shillor, \emph{Second order evolution
  equations with dynamic boundary conditions}, J. Math. Anal. Appl.
  \textbf{197} (1996), no.~3, 781--795.

\bibitem{barb}
S.~Barb, \emph{Topics in {G}eometric {A}nalysis with {A}pplications to
  {P}artial {D}ifferential {E}quations}, ProQuest LLC, Ann Arbor, MI, 2009,
  Thesis (Ph.D.)--University of Missouri - Columbia.

\bibitem{BergerMarcel}
M.~Berger, \emph{Geometry {I}}, Universitext, Springer-Verlag, Berlin, 2009,
  Translated from the 1977 French original by M. Cole and S. Levy, Fourth
  printing of the 1987 English translation [MR0882541].

\bibitem{Boothby}
W.~M. Boothby, \emph{An introduction to differentiable manifolds and
  {R}iemannian geometry}, Academic Press [A subsidiary of Harcourt Brace
  Jovanovich, Publishers], New York-London, 1975, Pure and Applied Mathematics,
  No. 63.

\bibitem{brezis2}
H.~Brezis, \emph{Functional analysis, {S}obolev spaces and partial differential
  equations}, Universitext, Springer, New York, 2011.

\bibitem{CalatroniColli}
L.~Calatroni and P.~Colli, \emph{Global solution to the {A}llen-{C}ahn equation
  with singular potentials and dynamic boundary conditions}, Nonlinear Anal.
  \textbf{79} (2013), 12--27.

\bibitem{CavalcantiKM}
M.~M. Cavalcanti and M.~Khemmoudj, A.and~Medjden, \emph{Uniform stabilization
  of the damped {C}auchy-{V}entcel problem with variable coefficients and
  dynamic boundary conditions}, J. Math. Anal. Appl. \textbf{328} (2007),
  no.~2, 900--930.

\bibitem{ColliFukao}
P.~Colli and T.~Fukao, \emph{Cahn-{H}illiard equation with dynamic boundary
  conditions and mass constraint on the boundary}, J. Math. Anal. Appl.
  \textbf{429} (2015), no.~2, 1190--1213.

\bibitem{conradmorgul}
F.~Conrad and {\"O}.~Morg{\"u}l, \emph{On the stabilization of a flexible beam
  with a tip mass}, SIAM J. Control Optim. \textbf{36} (1998), no.~6,
  1962--1986 (electronic).

\bibitem{Dambrine}
M.~Dambrine, D.~Kateb, and J.~Lamboley, \emph{An extremal eigenvalue problem
  for the {W}entzell-{L}aplace operator}, Ann. Inst. H. Poincar\'{e} C Anal.
  Non Lin\'{e}aire \textbf{33} (2016), no.~2, 409--450.

\bibitem{darmvanhorssen}
Darmawijoyo and W.~T. van Horssen, \emph{On boundary damping for a weakly
  nonlinear wave equation}, Nonlinear Dynam. \textbf{30} (2002), no.~2,
  179--191.

\bibitem{dautraylionsvol2}
R.~Dautray and J.-L. Lions, \emph{Mathematical analysis and numerical methods
  for science and technology. {V}ol. 2}, Springer-Verlag, Berlin, 1988,
  Functional and Variational Methods.

\bibitem{dautraylionsvol3}
\bysame, \emph{Mathematical analysis and numerical methods for science and
  technology. {V}ol. 3}, Springer-Verlag, Berlin, 1990, Spectral Theory and
  Applications.

\bibitem{doroninlarkinsouza}
G.~G. Doronin, N.~A. Lar{\cprime}kin, and A.~J. Souza, \emph{A hyperbolic
  problem with nonlinear second-order boundary damping}, Electron. J.
  Differential Equations (1998), No. 28, 10 pp. (electronic).

\bibitem{DuWangXia}
F.~Du, Q.~Wang, and C.~Xia, \emph{Estimates for eigenvalues of the
  {W}entzell-{L}aplace operator}, J. Geom. Phys. \textbf{129} (2018), 25--33.

\bibitem{Evans}
L.~C. Evans, \emph{Partial differential equations}, second ed., Graduate
  Studies in Mathematics, vol.~19, American Mathematical Society, Providence,
  RI, 2010.

\bibitem{FGGGR}
A.~Favini, C.~G. Gal, G.~Ruiz~Goldstein, J.~A. Goldstein, and S.~Romanelli,
  \emph{The non-autonomous wave equation with general {W}entzell boundary
  conditions}, Proc. Roy. Soc. Edinburgh Sect. A \textbf{135} (2005), no.~2,
  317--329.

\bibitem{gilbarg}
D.~Gilbarg and N.~Trudinger, \emph{Elliptic partial differential equations of
  second order}, 2nd ed., Springer--Verlag, Berlin--New York, 1983.

\bibitem{goldsteingisele}
G.~Ruiz Goldstein, \emph{Derivation and physical interpretation of general
  boundary conditions}, Adv. Differential Equations \textbf{11} (2006), no.~4,
  457--480.

\bibitem{grisvard}
P.~Grisvard, \emph{Elliptic problems in nonsmooth domains}, Classics in Applied
  Mathematics, vol.~69, Society for Industrial and Applied Mathematics (SIAM),
  Philadelphia, PA, 2011, Reprint of the 1985 original, with a foreword by
  Susanne C. Brenner.

\bibitem{guoxu}
B.~Guo and C.-Z. Xu, \emph{On the spectrum-determined growth condition of a
  vibration cable with a tip mass}, IEEE Trans. Automat. Control \textbf{45}
  (2000), no.~1, 89--93.

\bibitem{hebey}
E.~Hebey, \emph{Nonlinear analysis on manifolds: {S}obolev spaces and
  inequalities}, Courant Lecture Notes in Mathematics, vol.~5, New York
  University, Courant Institute of Mathematical Sciences, New York; American
  Mathematical Society, Providence, RI, 1999.

\bibitem{heminna1}
A.~Heminna, \emph{Stabilisation fronti\`ere de probl\`emes de {V}entcel}, C. R.
  Acad. Sci. Paris S\'{e}r. I Math. \textbf{328} (1999), no.~12, 1171--1174.

\bibitem{Heminna2}
\bysame, \emph{Stabilisation fronti\`ere de probl\`emes de {V}entcel}, ESAIM
  Control Optim. Calc. Var. \textbf{5} (2000), 591--622.

\bibitem{jost}
J.~Jost, \emph{Riemannian {G}eometry and {G}eometric {A}nalysis}, fifth ed.,
  Universitext, Springer-Verlag, Berlin, 2008.

\bibitem{KnopfLiu}
P.~Knopf and C.~Liu, \emph{On second-order and fourth-order elliptic systems
  consisting of bulk and surface {PDE}s: well-posedness, regularity theory and
  eigenvalue problems}, Interfaces Free Bound. \textbf{23} (2021), no.~4,
  507--533.

\bibitem{LeoniSobolev2}
G.~Leoni, \emph{A first course in {S}obolev spaces}, second ed., Graduate
  Studies in Mathematics, vol. 181, American Mathematical Society, Providence,
  RI, 2017.

\bibitem{LewickaPeres}
M.~Lewicka and Y.~Peres, \emph{Which domains have two-sided supporting unit
  spheres at every boundary point?}, Expo. Math. \textbf{38} (2020), no.~4,
  548--558.

\bibitem{LLXiao1}
C.~Li, J.~Liang, and T.-J. Xiao, \emph{Asymptotic behaviours of solutions for
  wave equations with damped {W}entzell boundary conditions but no interior
  damping}, J. Differential Equations \textbf{271} (2021), 76--106.

\bibitem{LLXiao2}
\bysame, \emph{Regularity and stability of wave equations with variable
  coefficients and {W}entzell type boundary conditions}, J. Differential
  Equations \textbf{374} (2023), 548--592.

\bibitem{LiXiao}
C.~Li and T.-J. Xiao, \emph{Asymptotics for wave equations with {W}entzell
  boundary conditions and boundary damping}, Semigroup Forum \textbf{94}
  (2017), no.~3, 520--531.

\bibitem{morgulraoconrad}
{\"O}.~Morg{\"u}l, B.~Rao, and F.~Conrad, \emph{On the stabilization of a cable
  with a tip mass}, IEEE Trans. Automat. Control \textbf{39} (1994), no.~10,
  2140--2145.

\bibitem{MugnoloCosine}
D.~Mugnolo, \emph{Operator matrices as generators of cosine operator
  functions}, Integral Equations Operator Theory \textbf{54} (2006), no.~3,
  441--464.

\bibitem{mugnolo2011}
\bysame, \emph{Damped wave equations with dynamic boundary conditions}, J.
  Appl. Anal. \textbf{17} (2011), no.~2, 241--275.

\bibitem{mugnvit}
D.~Mugnolo and E.~Vitillaro, \emph{The wave equation with acoustic boundary
  conditions on non-locally reacting surfaces}, to appear on Mem. Amer. Math.
  Soc.

\bibitem{NicaiseLaoubi}
S.~Nicaise and K.~Laoubi, \emph{Polynomial stabilization of the wave equation
  with {V}entcel's boundary conditions}, Math. Nachr. \textbf{283} (2010),
  no.~10, 1428--1438.

\bibitem{pucvit}
P.~Pucci and E.~Vitillaro, \emph{Approximation by regular functions in
  {S}obolev spaces arising from doubly elliptic problems}, Boll. Unione Mat.
  Ital. \textbf{13} (2020), no.~4, 487--494. \MR{4172949}

\bibitem{roman}
S.~Roman, \emph{Advanced linear algebra}, third ed., Graduate Texts in
  Mathematics, vol. 135, Springer, New York, 2008.

\bibitem{SprekelsWu}
J.~Sprekels and H.~Wu, \emph{A note on parabolic equation with nonlinear
  dynamical boundary condition}, Nonlinear Anal. \textbf{72} (2010), no.~6,
  3028--3048.

\bibitem{sternberg}
S.~Sternberg, \emph{Lectures on differential geometry}, second ed., Chelsea
  Publishing Co., New York, 1983, With an appendix by Sternberg and Victor W.
  Guillemin.

\bibitem{taylor}
M.~E. Taylor, \emph{Partial differential equations}, Texts in Applied
  Mathematics, vol.~23, Springer-Verlag, New York, 1996, Basic theory.

\bibitem{taylor3}
\bysame, \emph{Partial differential equations. {III}}, Applied Mathematical
  Sciences, vol. 117, Springer-Verlag, New York, 1997.

\bibitem{triebel}
H.~Triebel, \emph{Interpolation theory, function spaces, differential
  operators}, North-Holland, Amsterdam, 1978.

\bibitem{vazvitHLB}
J.~L. V\'{a}zquez and E.~Vitillaro, \emph{Heat equation with dynamical boundary
  conditions of reactive-diffusive type}, J. Differential Equations
  \textbf{250} (2011), no.~4, 2143--2161.

\bibitem{AMS}
E.~Vitillaro, \emph{Strong solutions for the wave equation with a kinetic
  boundary condition}, Recent trends in nonlinear partial differential
  equations. {I}. {E}volution problems, Contemp. Math., vol. 594, Amer. Math.
  Soc., Providence, RI, 2013, pp.~295--307.

\bibitem{Dresda1}
\bysame, \emph{On the {W}ave {E}quation with {H}yperbolic {D}ynamical
  {B}oundary {C}onditions, {I}nterior and {B}oundary {D}amping and {S}ource},
  Arch. Ration. Mech. Anal. \textbf{223} (2017), no.~3, 1183--1237.

\bibitem{Dresda2}
\bysame, \emph{On the wave equation with hyperbolic dynamical boundary
  conditions, interior and boundary damping and supercritical sources}, J.
  Differential Equations \textbf{265} (2018), no.~10, 4873--4941.

\bibitem{Dresda3}
\bysame, \emph{Blow-up for the wave equation with hyperbolic dynamical boundary
  conditions, interior and boundary nonlinear damping and sources}, Discrete
  Contin. Dyn. Syst. Ser. S \textbf{14} (2021), no.~12, 4575--4608.

\bibitem{MP}
\bysame, \emph{Nontrivial solutions for the {L}aplace equation with a nonlinear
  {G}oldstein-{W}entzell boundary condition}, Commun. Anal. Mech. \textbf{15}
  (2023), no.~4, 811--830.

\bibitem{xiaoliang2}
T.-J. Xiao and L.~Jin, \emph{Complete second order differential equations in
  {B}anach spaces with dynamic boundary conditions}, J. Differential Equations
  \textbf{200} (2004), no.~1, 105--136.

\bibitem{xiaoliang1}
T.-J. Xiao and J.~Liang, \emph{Second order parabolic equations in {B}anach
  spaces with dynamic boundary conditions}, Trans. Amer. Math. Soc.
  \textbf{356} (2004), no.~12, 4787--4809.

\bibitem{zahn}
J.~Zahn, \emph{Generalized {W}entzell boundary conditions and quantum field
  theory}, Ann. Henri Poincar\'{e} \textbf{19} (2018), no.~1, 163--187.

\end{thebibliography}
\def\cprime{$'$}
\providecommand{\bysame}{\leavevmode\hbox to3em{\hrulefill}\thinspace}
\providecommand{\MR}{\relax\ifhmode\unskip\space\fi MR }
\providecommand{\MRhref}[2]{%
  \href{http://www.ams.org/mathscinet-getitem?mr=#1}{#2}
}
\providecommand{\href}[2]{#2}

\end{document}